\documentclass[preprint,12pt]{elsarticle}

\usepackage{amsmath, amsthm, amssymb}
\usepackage{hyperref}
\usepackage[utf8]{inputenc}
\usepackage[T1]{fontenc}

\usepackage{xcolor}

\newtheorem{theorem}{Theorem}[section]
\newtheorem{proposition}[theorem]{Proposition}

\newcommand{\R}{\mathbb R}

\newcommand{\pv}{\operatorname{P.V}}

\newcommand{\M}{\mathbf{M}}
\newcommand{\bv}{\mathbf{v}}
\newcommand{\fvect}[1]{\mathbf{#1}}

\newtheorem{lemma}[theorem]{Lemma}
\newtheorem{corollary}[theorem]{Corollary}
\newtheorem{definition}{Definition}

\journal{Journal of Differential Equations}

\begin{document}

\begin{frontmatter}

\title{Counterexample to the second eigenfunction having one zero for a non-local Schr\"{o}dinger operator}

\author[1]{Ben Andrews}
\ead{ben.andrews@anu.edu.au}
\author[2, fn1]{Sophie Chen}
\ead{sophie.c.chen@gmail.com}

\address[1]{Mathematical Sciences Institute, Australian National University, Canberra, Australia}
\address[2]{Mathematical Sciences Institute, Australian National University, Canberra, Australia}

\fntext[fn1]{Independent Researcher, ACT, Australia}

\begin{abstract}
We demonstrate that the second eigenfunction of a perturbed fractional Laplace operator on a bounded interval can exhibit two sign changes, in stark contrast with the classical expectation that it should have exactly one zero. Our construction employs the analytic perturbation theory of Kato and Rellich for degenerate eigenvalues to analyse an infinite potential well eigenvalue problem, and then uses a compactness and energy‐minimisation argument to extend this counterexample to finite potential wells. Although our detailed analysis focuses on the case s = 1/2 (the Cauchy process), our approach indicates that similar phenomena occur for other rational values of s in (0, 1). At the time of writing, this result provides one of the first rigorous insights into the qualitative behaviour of eigenfunctions for perturbed nonlocal Schrödinger operators.
\end{abstract}


\end{frontmatter}

\section{Introduction and counterexample outline}\label{ch_4_background}

Let $\Omega$ be an open bounded set in $\R^n$ with a Lipschitz boundary. For $\epsilon>0$ and $s$ in $(0,1)$, consider the eigenvalue equation of a perturbed fractional Laplace operator with Dirichlet boundary data:
\begin{equation}\label{evalue_problem_perturbed_frac_Lap}
 \left\{\begin{aligned}
        ((-\Delta)^s+\epsilon V(x))u(x)&=\lambda u(x), &x\in\Omega\\
        u(x)&=0, &x\in\R^n\setminus \Omega,
       \end{aligned}
 \right.
\end{equation}
where $V:\Omega\to\R$ is a bounded potential and $(-\Delta)^s$ is the fractional Laplacian defined pointwise by
\[
(-\Delta)^su(x):=c_{n,s}\lim_{\epsilon\to 0^+}\int_{\R^n\setminus B_\epsilon(x)}\frac{(u(x)-u(y))}{\left|x-y\right|^{n+2s}}dy
\]
for smooth and compactly supported functions $u$. It is known that problem (\ref{evalue_problem_perturbed_frac_Lap}) has a countable set of eigenvalues that can be written as an increasing sequence $-\infty<\lambda_1<\lambda_2\leq\lambda_3\leq\dots\leq\lambda_k\leq\dots$ in the index $k$, and the corresponding family of eigenfunctions 
forms an orthonormal basis for $L^2(\Omega)$. The first eigenvalue is simple and has an eigenfunction that can be chosen to be strictly positive in $\Omega$, 
repeated eigenvalues have finite multiplicity, and each eigenvalue has a variational structure. Moreover, if the domain $\Omega$ also satisfies the exterior ball condition, i.e., there exists a positive radius $r>0$ such that each point $p$ in the boundary $\partial\Omega$ can be touched by a ball of radius $r$ in $\R^n\setminus\Omega$ \cite{ROS},  then the eigenfunctions $\{u_k\}_{k=1}^\infty$ are H\"{o}lder continuous on $\R^n$. A proof of the above facts can be found in \cite{Zh} and references contained therein.

Now let $\Omega$ be the line interval $I=(-1,1)$ so that problem (\ref{evalue_problem_perturbed_frac_Lap}) becomes
\begin{equation}\label{evalue_problem_perturbed_frac_Lap_1d}
 \left\{\begin{aligned}
        ((-\Delta)^s+\epsilon V(x))u(x)&=\lambda u(x), &x\in I\\
        u(x)&=0, &x\in\R\setminus I.
       \end{aligned}
 \right.
\end{equation}
It was shown by Ban\~{u}elos and Kulczycki in \cite{BK} that for $s=1/2$, the so-called `Cauchy process', the unperturbed eigenvalue problem
\begin{equation}\label{evalue_problem_unperturbed_frac_Lap_1d}
 \left\{\begin{aligned}
        (-\Delta)^su(x)&=\lambda u(x), &x\in I\\
        u(x)&=0, &x\in\R\setminus I
       \end{aligned}
 \right.
\end{equation}
has a second eigenvalue that is simple and the corresponding eigenfunction (modulo a constant) has exactly one zero. Furthermore, the second eigenfunction is antisymmetric, convex on one half of the interval, and concave on the other half. While these properties can be easily confirmed by direct computation for the classical Laplace operator with Dirichlet boundary condition on the interval (viz., the $s\to 1^-$ case), it is markedly more difficult to study the geometric properties of eigenfunctions of general fractional Dirichlet Laplace operators, even on the simplest geometry of a bounded line interval. Indeed, at the time of writing, many basic questions concerning the finer spectral properties of fractional Laplace operators remain unknown. 

Our aim in this paper is to construct a bounded potential on $I$ for which
the second eigenfunction of the Dirichlet fractional Laplace operator with
$s=1/2$ has at least two sign changes. To our knowledge, this is the first
such counterexample for a perturbed fractional Laplace operator. It demonstrates that the one-sign-change property of the second eigenfunction of the unperturbed problem is not, in general, preserved under bounded perturbations. 


\bigskip

\noindent\textbf{Counterexample outline}
\begin{figure}
    \centering
    \includegraphics[width=11cm]{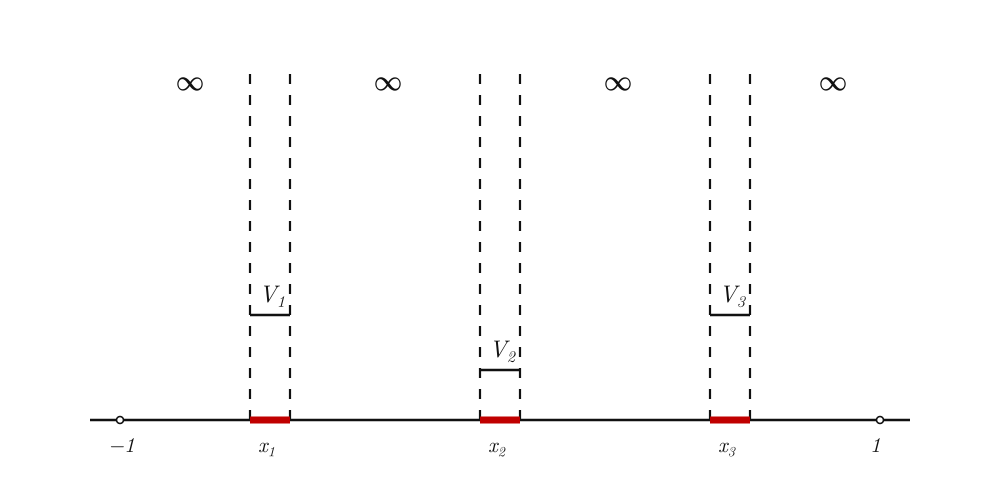}
    \caption{An example of the infinite potential well $V^\infty$. Red intervals represent $(x_i-\epsilon,x_i+\epsilon)$ for $i=1,2,3$ and small $\epsilon>0$.}
    \label{fig_infinite_potential}
\end{figure}

\bigskip

Before delving into the details of the counterexample in the remainder of this paper, we provide an overview of the three key steps to its construction.
\bigskip

\textbf{Step 1: Matrix eigenvalue problem.} First, consider an extreme case where the potential is bounded and constant only on a finite number of mutually disjoint, infinitesimally small sub-intervals $I_i=(x_i-\epsilon,x_i+\epsilon)$. 
The subintervals $I_i$ are centered at the ordered points $x_1<x_2<\ldots<x_k$ in $(-1,1)$ and are separated by distances that are large compared to their common length $2\epsilon$:
\[
\epsilon\ll \min\{|1-x_i|,|1+x_i|,|x_i-x_j|:j\neq i,i=1,\ldots,k\}.
\]
Off the sub-intervals the potential is infinite. Thus the singular potential of interest, which we henceforth refer to as the `infinite potential well', is 
\begin{equation}\label{potential_infinite}
V^{\infty}(x)=\left\{
\begin{array}{ll}
V_i, &x\in I_i,\qquad i=1,2,\ldots, k\\
\infty, &\mbox{else}.
\end{array}\right.
\end{equation}
See Figure \ref{fig_infinite_potential} for an example. To ensure that the energy associated with problem (\ref{evalue_problem_perturbed_frac_Lap_1d}) is finite when $V=V^\infty$, i.e.,
\[
\int_{\R}u(x)\left((-\Delta)^s+\epsilon V^\infty(x)\right)u(x)dx<\infty,
\]
the solution $u$ must be identically zero off the sub-intervals $I_i$ since $V^\infty$ is infinite there. Applying this fact to problem (\ref{evalue_problem_perturbed_frac_Lap_1d}) simplifies it to
\begin{equation}\label{simple_evalue_problem_perturbed_frac_Lap_1d}
 \left\{\begin{aligned}
        ((-\Delta)^s+\epsilon V_i)u(x)&=\lambda u(x), &x\in I_i,\qquad i=1,2,\ldots,k\\
        u(x)&=0, &x\in \R\setminus \bigcup_{i=1}^k I_i.
       \end{aligned}
 \right.
\end{equation}
As each subinterval $I_i$ is infinitesimally small, any admissible function $u$ can be approximated by a weighted sum of indicator functions on the $k$ subintervals:
\begin{align}\label{eigenvector_approximation}
    u(x)\approx \sum_{i=1}^k m_i\mathbf{1}_{I_i}(x),
\end{align}
where $m_i$ is a positive constant weight on the $i$th subinterval, and the indicator function $\mathbf{1_{A}}:\R\to\{0,1\}$ on the set $A\subset\R$ takes value one on $A$ but is otherwise zero. This means the perturbed operator in problem (\ref{simple_evalue_problem_perturbed_frac_Lap_1d}) has an approximate representation as a $k$-by-$k$ symmetric matrix $\M$ (coming from the operator being self-adjoint as we will later see) possessing the following properties: 
\begin{itemize}
    \item Diagonal elements are proportional to the finite potential values $V_i$ plus a term $\sigma_i$ related to the $L^2$-norm of the indicator function on $I_i$.
    \item Off-diagonal terms contain `interactions' $\sigma_i\sigma_j$, for $j\neq i$, related to the $L^2$-inner product between two indicator functions on distinct subintervals, and they can be assumed to be negative and decreasing in magnitude as we move away from the diagonal in both the vertical and horizontal directions.
\end{itemize}
So the matrix $\M$ has the form\\
\begin{equation}\label{matrix_putative_counterexample}
\M=\begin{bmatrix}
V_1+\sigma_1        & \sigma_1\sigma_2  & \dots             & \dots                 & \sigma_1\sigma_k      \\
\sigma_2\sigma_1    & \ddots            &                   &                       & \vdots                \\
\vdots              &                   & \ddots            &                       & \vdots                \\
\vdots              &                   &                   & \ddots                & \sigma_{k-1}\sigma_k  \\
\sigma_k\sigma_1    & \dots             & \dots             & \sigma_k\sigma_{k-1}  & V_k+\sigma_k
\end{bmatrix},
\end{equation}
with $\sigma_i\sigma_j=\sigma_j\sigma_i$ for all $j\neq i$ and $i=1,\ldots,k$. The approximate eigenfunctions are then given by the $k$-dimensional eigenvectors of $\M$. In the simplest case, where $k=3$, corresponding to $\M$ being a three-by-three matrix with eigenvectors that can have at most two sign changes, the problem can be fully analysed. We do this in Section \ref{putative_counterexample}. By studying this finite-dimensional matrix eigenvalue problem, we find that the second eigenvector has two sign changes provided the potential values $V_i$ is non-convex in $i$, i.e., $V_2<V_1,V_3$.\\

\textbf{Step 2: Infinite potential well counterexample.} Next, we produce a counterexample to the second eigenfunction of the infinite potential well eigenvalue problem (\ref{simple_evalue_problem_perturbed_frac_Lap_1d}) having only one sign change. The way to this result is to analyse the latter eigenvalue problem as a regular perturbation problem. This ultimately allows us to reduce it to the matrix eigenvalue problem solved in Step 1, so leading to a counterexample.\\

\textbf{Step 3: Finite potential well counterexample.} 
Finally, let each of the disjoint subintervals $I_i$ in $I$ have fixed uniform widths that are small compared to their distances apart. Take a sequence of finite deep potential wells that converge pointwise to the infinite potential $V^\infty$ of Step 2 as the well depth grows without bound.
We then proceed to produce a counterexample to the second eigenfunction of eigenvalue problem (\ref{evalue_problem_perturbed_frac_Lap}) having only one sign change, in the case where the bounded potential is a sufficiently deep finite potential well, using a compactness and energy minimisation argument.\\

Preliminary notions are defined in the next section. 

\section{Functional analytic framework}\label{subs_func_an}
For any $s$ in $(0,1)$, denote by $H^s(\R^n)$ the standard fractional Sobolev space
\begin{align}\label{def_std_frac_Sob_sp}
    H^s(\R^n)=\left\{u\in L^2(\R^n): [u]^2_{H^s(\R^n)}<\infty\right\}.
\end{align}
This is a Banach space equipped with the norm 
\begin{align}\label{def_H_s_norm}
    \|u\|_{H^s(\R^n)}=\left(\|u\|^2_{L^2(\R^n)}+\frac{c_{n,s}}{2}[u]^2_{H^s(\R^n)}\right)^{1/2},
\end{align}
where the function $u$ belongs to $H^s(\R^n)$, $\|u\|_{L^2(\R^n)}=\left(\int_{\R^n}|u|^2dx\right)^{1/2}$ is the $L^2$-norm, and the Gagliardo seminorm is given by
\begin{align}\label{def_Gagliardo_seminorm}
    [u]_{H^s(\R^n)}=\left(\int_{\R^n}\int_{\R^n}\frac{|u(x)-u(y)|^2}{|x-y|^{n+2s}}dxdy\right)^{1/2}.
\end{align}
\cite[p. 5]{NPV} The norm (\ref{def_H_s_norm}) is induced by the inner product
\begin{align}\label{def_H_s_IP}
    \langle u,v \rangle_{H^s(\R^n)}&=\langle u, v\rangle_{L^2(\R^n)}+\frac{c_{n,s}}{2}\iint\limits_{\R^n\times\R^n}\frac{(u(x)-u(y))(v(x)-v(y))}{|x-y|^{n+2s}}dxdy
\end{align}
for any functions $u,v$ in $H^s(\R^n)$. Thus $H^s(\R^n)$ with the real inner product (\ref{def_H_s_IP}) is a Hilbert space. 
 
Let $\Omega$ be an open and bounded set in $\R^n$. Define the subspace 
\begin{align}\label{def_subsp_H_0}
    H^s_0(\Omega):=\left\{u\in H^s(\R^n):u=0 \mbox{ a.e. on }\R^n\setminus\Omega\right\}
\end{align} 
and the symmetric bilinear form
\begin{align}\label{bilinear_form_res_frac_Lap_full}
\mathcal E_{\mathrm{res}}[u,v]
:=
\frac{c_{n,s}}{2}
\iint_{\mathbb R^n\times\mathbb R^n}
\frac{(u(x)-u(y))(v(x)-v(y))}
{|x-y|^{n+2s}}\,dx\,dy
\end{align}
for \(u,v\in H_0^s(\Omega)\). The zero exterior condition gives the equivalent representation
\begin{align}\label{bilinear_form_res_frac_Lap_zero}
\mathcal E_{\mathrm{res}}[u,v]
={}&
\frac{c_{n,s}}{2}
\iint_{\Omega\times\Omega}
\frac{(u(x)-u(y))(v(x)-v(y))}
{|x-y|^{n+2s}}\,dx\,dy \nonumber\\
&+
\int_\Omega \kappa(x)u(x)v(x)\,dx,
\end{align}
where
\begin{align}
\kappa(x)
=
c_{n,s}
\int_{\mathbb R^n\setminus\Omega}
\frac{1}{|x-y|^{n+2s}}\,dy,
\qquad x\in\Omega.
\end{align}
It is well known that $H^s_0(\Omega)$ with inner product
\begin{align*}
    \langle u,v\rangle_{H^s_0(\Omega)}:=\mathcal E_{\mathrm{res}}[u,v]
\end{align*}
is a Hilbert space with norm
\begin{align}\label{Hs0_norm}
    \|u\|_{H^s_0(\Omega)}=\mathcal{E}_{res}[u, u]^{\frac{1}{2}}=\left(\frac{c_{n,s}}{2}\right)^{\frac{1}{2}}[u]_{H^s(\R^n)}.
\end{align} 
By setting $v=u$ in (\ref{bilinear_form_res_frac_Lap_zero}), the quadratic form  $\mathcal{E}_{res}[u,u]$ is equivalently given by
\begin{align}\label{bilinear_form_res_frac_Lap}
\mathcal{E}_{res}[u,u]=\frac{c_{n,s}}{2}[u]^2_{H^s(\Omega)}+\int_{\Omega}\kappa(x)u(x)^2dx.
\end{align} \cite[Sec. 2.1]{BBDG}

On $H^s_0(\Omega)$, the norm (\ref{Hs0_norm}) is equivalent to the induced norm from $H^s(\R^n)$ specified by equation (\ref{def_H_s_norm}). The key step to showing this is the following lemma, which shows that $H^s_0(\Omega)$ is continuously embedded in $L^2(\Omega)$.

\begin{lemma}\label{lem_equiv_norms} 
    There exists a constant $c>0$ such that for any function $u$ in $H^s_0(\Omega)$,
    \[
        [u]_{H^s(\R^n)}\geq c\|u\|_{L^2(\Omega)}.
    \] 
\end{lemma}
\begin{proof} To generate a contradiction, suppose that for every constant $c>0$, there exists a function $u$ in $H^s_0(\Omega)$ such that 
\[
\|u\|_{L^2(\Omega)}>\frac{1}{c}[u]_{H^s(\R^n)}.
\]
In particular, there exists a sequence of constants $c_i$ going to zero and a sequence of functions $u_i$ in $H^s_0(\Omega)$, depending on $c_i$, such that $[u_i]_{H^s(\R^n)}<\infty$, but 
\[
\frac{1}{c_i}[u_i]_{H^s(\R^n)}<\|u_i\|_{L^2(\Omega)}\to\infty \qquad\mbox{as}\qquad i\to\infty\mbox{ (and }c_i\to 0).
\]
Equivalently,
\begin{align}\label{seminorm_limit}
        0\leq [u_i]_{H^s(\R^n)}<c_i\|u_i\|_{L^2(\Omega)}\to 0\qquad\mbox{as}\qquad i\to\infty\mbox{ (and }c_i\to 0).
\end{align}
We can therefore re-scale the $u_i$'s so that $[u_i]_{H^s(\R^n)}\to 0$ as $i\to\infty$, while $\|u_i\|_{L^2(\Omega)}=1$ for each $i$.

Now fix a sufficiently large $R>0$ such that $\Omega\subset B_R(0)\subset\R^n$, where $B_R(0)$ is the open ball of radius $R$ centred at zero in $\R^n$. This is an open Lipschitz domain with bounded boundary and so is a bounded extension domain for $H^s(B_R(0))$ by \cite[p. 33]{NPV}. Moreover, $\{u_i\}_{i\geq 1}$ is a bounded subset of $L^2(B_R(0))$, since $\|u\|_{L^2(B_R(0))}=\|u\|_{L^2(\Omega)}=1$ for each $i$. Finally,
\[
\sup_{u_i\in\{u_i\}_{i\geq 1}}[u_i]^2_{H^s(B_R(0))}=\sup_{u_i\in\{u_i\}_{i\geq 1}}\iint\limits_{B_R(0)\times B_R(0)}\frac{(u_i(x)-u_i(y))^2dxdy}{|x-y|^{n+2s}}<\infty
\]
as $[u_i]_{H^s(B_R(0))}\leq [u_i]_{H^s(\R^n)}\to 0$ as $i\to\infty$ from inequality (\ref{seminorm_limit}). It therefore follows from \cite[Th. 7.1]{NPV} that $\{u_i\}_{i\geq 1}$ is precompact in $L^2(B_R(0))$; i.e., there exists a subsequence $u_{i_j}$ that converges in $L^2$ to a limit point, $\bar{u}$, say, with $\|\bar{u}\|_{L^2(B_R(0))}=1$. On the other hand, inequality (\ref{seminorm_limit}) implies that $[\bar{u}]_{H^s(\R^n)}=0$. Explicitly, we have
\begin{align*}
    0=[\bar{u}]^2_{H^s(\R^n)}=[\bar{u}]^2_{H^s(B_R(0))}+2\int_{B_R(0)}\tilde{\kappa}(x)\bar{u}(x)^2dx,
\end{align*}
where $\tilde{\kappa}(x)=\int_{\R^n\setminus B_R(0)}\frac{1}{|x-y|^{n+2s}}dy$. The seminorm $[\bar{u}]_{H^s(B_R(0))}$ taking value zero implies that $\bar{u}$ is constant (a.e), while the final integral term being zero implies that the constant is zero. This contradicts the $L^2$ norm of $\bar{u}$ being one.
\end{proof}

\begin{proposition}\label{prop_equiv_norms} There exist positive constants $c$ and $C$ such that for any function $u$ in $H^s_0(\Omega)$,
\[
c\|u\|_{H^s_0(\Omega)}\leq \|u\|_{H^s(\R^n)}\leq C\|u\|_{H^s_0(\Omega)}.
\]
\end{proposition}
\begin{proof} The result follows from applying the respective definitions of $\|\cdot\|_{H^s_0(\Omega)}$ and $\|\cdot\|_{H^s(\R^n)}$ and Lemma \ref{lem_equiv_norms}.
\end{proof}

Using standard arguments, one can show that $H^s_0(\Omega)$ is both closed and dense in $L^2(\Omega)$.\\

The function space $H^s(\R^n)$ is related to the fractional Laplace operator $(-\Delta)^s$ on $\R^n$. Similarly, the subspace $H^s_0(\Omega)$ is associated with its restriction to $\Omega\subset\R^n$ with Dirichlet boundary condition on $\R^n\setminus\Omega$. This latter operator is called the `restricted fractional Laplace operator' and is denoted by $(-\Delta)^s_{res}$. As an $L^2$-valued operator, its domain is
\begin{align}\label{operator_domain_res_frac_lap}
D((-\Delta)^s_{\mathrm{res}}):=\big\{u\in H^s_0(\Omega):\;&\text{there exists }f\in L^2(\Omega)
\text{ such that} \nonumber\\
&\mathcal{E}_{\mathrm{res}}[u,v]
 =\langle f,v\rangle_{L^2(\Omega)}
 \text{ for every }v\in H^s_0(\Omega)\big\}.
\end{align}
For $u\in D((-\Delta)^s_{\mathrm{res}})$, the function $f$ is unique, as
$H^s_0(\Omega)$ is dense in $L^2(\Omega)$ and the $L^2$-inner product is continuous, and we define
$(-\Delta)^s_{res} u=f$. Thus
\begin{align}\label{operator_form_identity}
    \mathcal{E}_{\mathrm{res}}[u,v]
    =\langle (-\Delta)^s_{res} u,v\rangle_{L^2(\Omega)},
    \qquad u\in D((-\Delta)^s_{\mathrm{res}}),\quad v\in H^s_0(\Omega).
\end{align}
This domain is dense in $L^2(\Omega)$; see below.

The restricted fractional Laplacian is self-adjoint. 

\begin{lemma}\label{lemma_self_adj_res_frac_Lapl}
Let $\Omega$ be an open bounded subset of $\mathbb{R}^n$. The restricted
fractional Laplace operator
\[
    (-\Delta)^s_{\mathrm{res}}
    :
    D((-\Delta)^s_{\mathrm{res}})\subset L^2(\Omega)
    \longrightarrow L^2(\Omega),
\]
with operator domain defined by
\eqref{operator_domain_res_frac_lap}, is positive and self-adjoint.
\end{lemma}
\begin{proof} We first show that $(-\Delta)^s_{\mathrm{res}}$ is bijective. Let
$f\in L^2(\Omega)$ and define
\[
    F_f(v):=\langle f,v\rangle_{L^2(\Omega)},
    \qquad v\in H^s_0(\Omega).
\]
By the Cauchy--Schwarz inequality and
Lemma~\ref{lem_equiv_norms},
\[
    |F_f(v)|
    \leq
    \|f\|_{L^2(\Omega)}
    \|v\|_{L^2(\Omega)}
    \leq
    C\|f\|_{L^2(\Omega)}
    \|v\|_{H^s_0(\Omega)}.
\]
$F_f$ is therefore a bounded linear functional on the Hilbert space
$H^s_0(\Omega)$. By the Riesz representation theorem, there exists a
unique $u\in H^s_0(\Omega)$ such that
\[
    \mathcal E_{\mathrm{res}}[u,v]
    =
    \langle f,v\rangle_{L^2(\Omega)}
    \qquad
    \text{for every }v\in H^s_0(\Omega).
\]
Recalling the definition of $D((-\Delta)^s_{\mathrm{res}})$, this is equivalent to
\[
    u\in D((-\Delta)^s_{\mathrm{res}}),
    \qquad
    (-\Delta)^s_{\mathrm{res}}u=f.
\]
Existence and uniqueness for every $f\in L^2(\Omega)$ show
that $(-\Delta)^s_{\mathrm{res}}$ is bijective.

Next, the operator domain $D((-\Delta)^s_{\mathrm{res}})$ is dense in $L^2(\Omega)$. Indeed,
suppose that $h\in L^2(\Omega)$ is orthogonal to $D((-\Delta)^s_{\mathrm{res}})$.
Since $(-\Delta)^s_{\mathrm{res}}$ is onto, there exists
$u\in D((-\Delta)^s_{\mathrm{res}})$ such that
$(-\Delta)^s_{\mathrm{res}}u=h$. Hence,
\[
    0
    =
    \langle h,u\rangle_{L^2(\Omega)}
    =
    \langle (-\Delta)^s_{\mathrm{res}}u,u\rangle_{L^2(\Omega)}
    =
    \mathcal E_{\mathrm{res}}[u,u]
    =
    \|u\|_{H^s_0(\Omega)}^2.
\]
Consequently $u=0$, so
$h=(-\Delta)^s_{\mathrm{res}}u=0$. It follows that
$D((-\Delta)^s_{\mathrm{res}})$ is dense in $L^2(\Omega)$.

For any $u,v\in D((-\Delta)^s_{\mathrm{res}})$,
\begin{align*}
    \langle (-\Delta)^s_{\mathrm{res}}u,v\rangle_{L^2(\Omega)}
    &=
    \mathcal E_{\mathrm{res}}[u,v] \\
    &=
    \mathcal E_{\mathrm{res}}[v,u] \\
    &=
    \langle u,(-\Delta)^s_{\mathrm{res}}v\rangle_{L^2(\Omega)}.
\end{align*}
Thus $(-\Delta)^s_{\mathrm{res}}$ is symmetric. It is also positive,
since
\[
    \langle (-\Delta)^s_{\mathrm{res}}u,u\rangle_{L^2(\Omega)}
    =
    \mathcal E_{\mathrm{res}}[u,u]
    \geq 0.
\]

It remains to show that $(-\Delta)^s_{\mathrm{res}}$ is self-adjoint.
Since $D((-\Delta)^s_{\mathrm{res}})$ is dense in $L^2(\Omega)$ and symmetric, this implies that  
\[
    D((-\Delta)^s_{\mathrm{res}})
    \subset
    D\bigl({(-\Delta)^s_{\mathrm{res}}}^*\bigr)
\]
and the two operators agree on $D((-\Delta)^s_{\mathrm{res}})$. So we only need to prove the reverse inclusion of the
domains.

Let $v\in D\bigl({(-\Delta)^s_{\mathrm{res}}}^*\bigr)$. As $(-\Delta)^s_{\mathrm{res}}$ is onto $L^2(\Omega)$, there exists
$u\in D((-\Delta)^s_{\mathrm{res}})$ such that
\[
    (-\Delta)^s_{\mathrm{res}}u
    =
    {(-\Delta)^s_{\mathrm{res}}}^*v.
\]
Then, for every $x\in D((-\Delta)^s_{\mathrm{res}})$,
\begin{align*}
    \langle (-\Delta)^s_{\mathrm{res}}x,v\rangle_{L^2(\Omega)}
    &=
    \langle x,{(-\Delta)^s_{\mathrm{res}}}^*v\rangle_{L^2(\Omega)} \\
    &=
    \langle x,(-\Delta)^s_{\mathrm{res}}u\rangle_{L^2(\Omega)} \\
    &=
    \langle (-\Delta)^s_{\mathrm{res}}x,u\rangle_{L^2(\Omega)},
\end{align*}
where the first equality follows from the definition of the adjoint
and the last from symmetry. Hence,
\[
    \langle (-\Delta)^s_{\mathrm{res}}x,v-u\rangle_{L^2(\Omega)}
    =
    0
    \qquad
    \text{for every }x\in D((-\Delta)^s_{\mathrm{res}}).
\]
As $(-\Delta)^s_{\mathrm{res}}$ is onto, the vector
$(-\Delta)^s_{\mathrm{res}}x$ ranges over all of $L^2(\Omega)$.
Therefore, writing $f:=(-\Delta)^s_{\mathrm{res}}x$,
\[
    \langle f,v-u\rangle_{L^2(\Omega)}
    =
    0
    \qquad
    \text{for every }f\in L^2(\Omega).
\]
Taking $f=v-u$ yields
\[
    \|v-u\|_{L^2(\Omega)}^2=0,
\]
and so $v=u\in D((-\Delta)^s_{\mathrm{res}})$. We thus conclude that
\[
    D\bigl({(-\Delta)^s_{\mathrm{res}}}^*\bigr)
    \subset
    D((-\Delta)^s_{\mathrm{res}}).
\]

Together with the reverse inclusion obtained from symmetry,
\[
    D\bigl({(-\Delta)^s_{\mathrm{res}}}^*\bigr)
    =
    D((-\Delta)^s_{\mathrm{res}}),
\]
and the two operators agree on this common domain. Thus,
\[
    (-\Delta)^s_{\mathrm{res}}
    =
    {(-\Delta)^s_{\mathrm{res}}}^*,
\]
showing that $(-\Delta)^s_{\mathrm{res}}$ is self-adjoint.
\end{proof}

Absorb the non-local Dirichlet boundary condition of
\eqref{evalue_problem_perturbed_frac_Lap} into the restricted fractional
Laplace operator, and define
\[
    L_V:=(-\Delta)^s_{\mathrm{res}}+\epsilon M_V
    =
    (-\Delta)^s_{\mathrm{res}}+\epsilon M_V
\]
on $L^2(\Omega)$. Since $V\in L^\infty(\Omega)$ is real-valued,
multiplication by $V$ gives a bounded self-adjoint operator $M_V$ on
$L^2(\Omega)$. Hence
\[
    D(L_V)=D((-\Delta)^s_{\mathrm{res}}).
\]

A real number $\lambda$ is an eigenvalue of
\eqref{evalue_problem_perturbed_frac_Lap} if there exists a non-zero
function $u\in D((-\Delta)^s_{\mathrm{res}})$ such that
\[
    L_Vu=\lambda u
    \qquad\text{in }L^2(\Omega).
\]
The function $u$ is then an eigenfunction corresponding to $\lambda$. Using \eqref{operator_form_identity}, this operator eigenvalue equation is
equivalent to the variational formulation
\begin{equation}\label{weak_form_perturbed_frac_lap}
    \mathcal E_{\mathrm{res}}[u,w]
    +
    \epsilon\int_\Omega V(x)u(x)w(x)\,dx
    =
    \lambda\int_\Omega u(x)w(x)\,dx
\end{equation}
for every $w\in H^s_0(\Omega)$.\\

Equipped with the requisite functional analytic framework, we now proceed to the first step of constructing the counterexample

\section{Matrix eigenvalue problem}\label{putative_counterexample}
Following the strategy in Section \ref{ch_4_background}, we use the matrix $\M$ given by equation (\ref{matrix_putative_counterexample}) to approximate the infinite-dimensional eigenvalue problem (\ref{simple_evalue_problem_perturbed_frac_Lap_1d}). Our goal is to completely analyse the simplest such approximation, corresponding to $k=3$, and show that it has a second eigenfunction with two sign changes under certain conditions on the discrete potential $V_i$ on each of the subintervals $I_i$.

We can represent this three-by-three model using the simplified matrix
\begin{equation}\label{matrix_three-by-three}
    \bar{\M}=\begin{bmatrix}
    U & c & b\\
    c & V & a\\
    b & a & W
    \end{bmatrix},
\end{equation}
where $U, V, W$ denote the diagonal `potential' terms $V_i+\sigma_i$, and the off-diagonal elements $a,b,c$ are the interactions $\sigma_i\sigma_j$, with $j\neq i$, in $\M$. As before, we assume $a,b,c<0$ and $|a|, |c|>|b|$. 

For any vector $\bv$ in $\R^3$ with (Euclidean) norm $\|\bv\|=1$, we can define the Rayleigh-Ritz quotient $\mathcal{E}:\R^n\setminus\{0\}\to\R$ associated with $\bar{\M}$, the so-called `energy':
\begin{align}\label{energy_putative_counterexample}
\mathcal{E}(\bv)&=\bv^T\bar{\M}\bv,
\end{align}
where $\bv^T$ denotes the transpose of the column vector $\bv$. Notice that when $\bv$ is an eigenvector of $\bar{\M}$, $\mathcal{E}(\bv)$ returns its corresponding eigenvalue.

\begin{theorem}\label{positive_evector} The minimum eigenvalue $\lambda_M$ of the matrix $\bar{\M}$ is simple and has an eigenvector $\bv_M$ with all components positive.
\end{theorem}
\begin{proof} Let $\mathbf{A}=\sigma \mathbf{I}-\bar{\M}$, where $\sigma>0$ is a real constant and $\mathbf{I}$ denotes the three-by-three identity matrix. Observe that if $\lambda$ is an eigenvalue of $\bar{\M}$ with corresponding eigenvector $\bv$, then 
\[
\mathbf{A}\bv=(\sigma \mathbf{I}-\bar{\M})\bv=\sigma \mathbf{I}\bv-\bar{\M}\bv=\sigma \bv-\lambda \bv=(\sigma-\lambda)\bv.
\]
That is, $\bv$ is also an eigenvector of $\mathbf{A}$ but corresponding to the eigenvalue $\mu:=\sigma-\lambda$.

If we let $A_{ij}$ and $\bar M_{ij}$ denote the $ij$-th entry of $\mathbf{A}$ and and $\bar{\M}$, respectively, then $A_{ij}=\sigma \delta_{ij}-\bar{M}_{ij}$, where $\delta_{ij}$ is the Kronecker delta function. The diagonal entries $A_{ii}=\sigma-\bar{M}_{ii}$ are strictly bigger than zero if we choose $\sigma>\max_i \{M_{ii}\}$. Off-diagonal entries $A_{ij}=-\bar{M}_{ij}$ are positive for all $i\neq j$, since the off-diagonal entries of $\bar{\M}$ are negative by construction. Thus $\mathbf{A}$ is a strictly positive real square matrix provided that $\sigma>\max_{i}\{\bar{M}_{ii}\}$. 

By the Perron-Frobenius theorem, $\mathbf{A}$ has a largest eigenvalue $\rho$ that is simple and a corresponding eigenvector $\bv_{\rho}$ with all components positive. Moreover, there are no other positive eigenvectors other than positive multiples of $\bv_{\rho}$. Since the eigenvalues of $\bar{\M}$ are related to those of $\mathbf{A}$ by $\lambda=\sigma - \mu$, we see that if $\rho$ is the largest eigenvalue of $A$, then $\sigma-\rho$ gives the smallest eigenvalue of $\bar{\M}$. Set $\lambda_M:=\sigma-\rho$. As $\rho$ is simple, $\lambda_M$ is also simple. Furthermore, the eigenvector $\bv_\rho$ corresponding to $\rho$ is the same eigenvector for $\lambda_M$. So $\bv_M:=\bv_\rho$ has all components positive.
\end{proof}

Having established that the smallest eigenvalue $\lambda_M$ of $\bar{\M}$ is simple, a corollary is that repeated eigenvalues occur only if the remaining two eigenvalues $\lambda_2$ and $\lambda_3$ are equal. As the $\lambda_M$-eigenspace is a straight-line through the origin along a unit vector in the positive octant of $\R^3$, its orthogonal complement $\bv_M^\perp$ is a plane through the origin in $\R^3$. The spectral theorem for real symmetric matrices guarantees that we can find an orthonormal set of eigenvectors $\{\bv_M,\bv_2,\bv_3,\}$ of $\bar \M$ that form a basis for $\R^3$. So $\bv_M^\perp$ is spanned by $\bv_2,\bv_3$ and contains any eigenvector with a zero component.

\begin{corollary}\label{second_eigenvector}
Any multiplicity of eigenvalues occurs if $\lambda_2=\lambda_3$, and $\bv_M^\perp=\operatorname{span}\{\bv_2, \bv_3\}$ contains eigenvectors of the form $[x,y,0], [0,y,z]$, and $[x,0,z]$.
\end{corollary}

We can solve the eigenvalue equation
\begin{align}\label{evalue_equation}
\bar{\M}\bv=\lambda \bv
\end{align}
for eigenvectors having the forms shown in Corollary \ref{second_eigenvector} corresponding to an eigenvalue $\lambda$ representing $\lambda_2$ or $\lambda_3$. If $\bv$ has the form $[x,y,0]$, a small computation gives
\begin{align}\label{evector_calc}
    0&  =(\bar{\M}-\lambda \mathbf{I})\begin{bmatrix}
           x \\
           y \\
           0
         \end{bmatrix}\nonumber\\
        &=  \begin{bmatrix}
            U-\lambda & c & b\\
            c & V-\lambda & a\\
            b & a & W-\lambda
            \end{bmatrix}
            \begin{bmatrix}
                x \\
                y \\
                0
            \end{bmatrix}\nonumber\\
        &=  \begin{bmatrix}
                (U-\lambda)x+cy\\
                cx +(V-\lambda)y\\
                bx+ay
            \end{bmatrix}.
\end{align}
From the last entry in equation (\ref{evector_calc}), we can choose $x=a$ and $y=-b$, so that 
\begin{align}\label{evector_solution}
\bv=\begin{bmatrix}
    x\\
    y\\
    0
\end{bmatrix} = \begin{bmatrix}
                    a\\
                    -b\\
                    0
                \end{bmatrix}, a^2+b^2=1;
\end{align}
i.e., eigenvectors of the form $[x,y,0]$ are scalar multiplies of $[a,-b,0]$. 
The potential values $U$ and $V$ at which the eigenvector $[a,-b,0]$ occurs are
\begin{align}\label{potential_U}
    U&=\lambda+\frac{bc}{a}
\end{align}
and
\begin{align}\label{potential_V}
    V&=\lambda+\frac{ac}{b}.
\end{align}
A similar calculation shows that eigenvectors of the form $[0,y,z]$ are multiples of $[0,-b,c]$ and occur when 
\begin{align}\label{potential_W}
V=\lambda+\frac{ac}{b} \qquad\mbox{ and }\qquad W=\lambda+\frac{ab}{c}.
\end{align}
While eigenvectors of the form $[x,0,z]$ are multiples of $[a,0,-c]$ and occur when $U=\lambda+\frac{bc}{a}$ and $W=\lambda+\frac{ab}{c}$. Notice that $[a,0,-c]=[a,-b,0]-[0,-b,c]$, so we have two linearly independent eigenvectors, $[a,-b,0]$ and $[0,-b,c]$, corresponding to $\lambda$. This means we can narrow our focus to the latter two linearly independent eigenvectors for the remainder of the calculations.\\

First we check the degeneracy of $\lambda$: Under what conditions is $\lambda$ degenerate, $\lambda_2=\lambda_3$, and when is it simple, $\lambda_2\neq \lambda_3$? To answer this question, we make the following key observation: adding a scalar multiple of the identity matrix $\beta \mathbf{I}$ to $\bar{\M}$ implies that if $\bar{\M}+\beta \mathbf{I}$ has an eigenvalue $\mu$, then $\bar{\M}$ has an eigenvalue $\lambda=\mu-\beta$. Choosing $\beta=\mu$ gives zero as an eigenvalue of $\bar{\M}$. Note that adding a scalar multiple of $\mathbf{I}$ to $\bar{\M}$ does not change the eigenvectors.

\begin{figure}[htp]
    \centering
    \includegraphics[width=12cm]{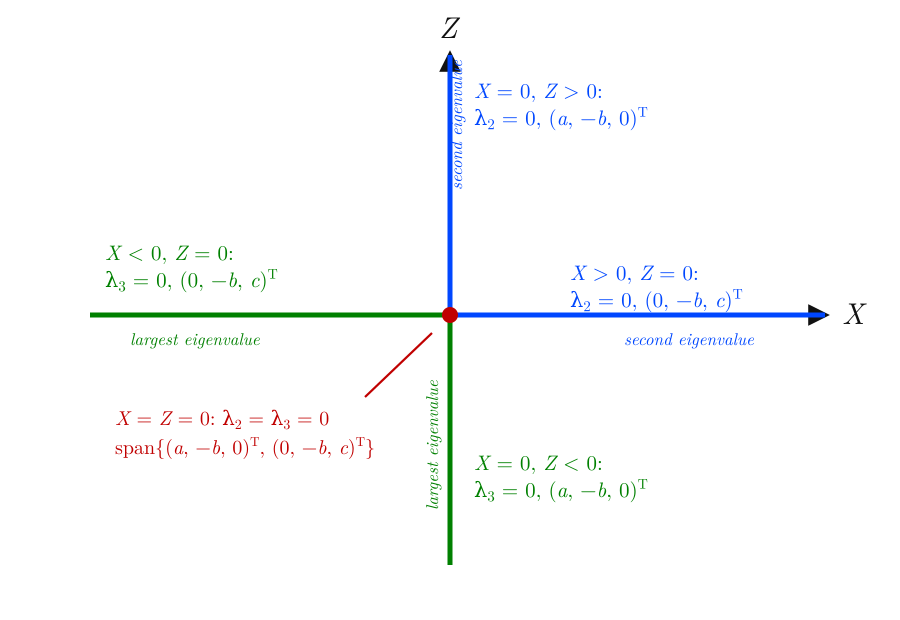}
    \caption{The second eigenvalue $\lambda$ is degenerate only at the origin. It is simple with eigenvalue $\lambda=0$ on the $Z$ and $X$ axes, respectively, excluding the origin.}\label{lambda_degeneracy}
\end{figure}

By adding a multiple of $\mathbf{I}$ to $\bar{\M}$, so that zero is an eigenvalue of $\bar{\M}$ and
\begin{align}\label{transformed_V}
V=\frac{ac}{b},
\end{align}
and setting 
\begin{align}\label{transformed_U_W}
    U=\frac{bc}{a}+X,\qquad W=\frac{ab}{c}+Z,
\end{align}
where $X$ and $Z$ are free variables, we end up with the following cases:
\begin{itemize}
    \item $\lambda$ is degenerate only when $X=Z=0$, and this gives $\lambda_2=\lambda_3=0$, with corresponding linearly independent eigenvectors $[a,-b,0]$ and $[0,-b,c]$.
    \item $\lambda$ is non-degenerate when 
        \begin{itemize}
            \item $X=0$ and $Z\neq 0$, and then $\lambda=0$ with corresponding eigenvector $[a,-b,0]$; or
            \item $Z=0$ and $X\neq 0$, and $\lambda=0$ with corresponding eigenvector $[0,-b,c]$.
        \end{itemize}
\end{itemize}
The above conclusions result from comparing the original expressions for each of the potential values $U, V,W$ in equations (\ref{potential_U})-(\ref{potential_W}) to their transformed values in equations (\ref{transformed_V}) and (\ref{transformed_U_W}). See Figure \ref{lambda_degeneracy} for a visualisation of these possibilities.\\

Next, we need to work out the conditions under which each of the eigenvectors $[0,-b,c]$ and $[a,-b,0]$ corresponds to $\lambda=0$ as the second (non-degenerate) eigenvalue of $\bar{\M}$. Since the computation is analogous in both cases, we only provide details for one of the cases. 

Take the eigenvector $[0,-b,c]$ corresponding to $\lambda=0$ on $\{Z=0\}\cap\{X\neq 0\}:=\{(X,0): X\in\R\setminus\{0\}\}$, with the finite potentials taking values $V=\frac{ac}{b}, U=\frac{bc}{a}+X, \mbox{ and } W=\frac{ab}{c}$ in the matrix $\bar{\M}$ given in equation (\ref{matrix_three-by-three}): 
\begin{equation}\label{matrix_three-by-three_lambda_zero}
    \bar{\M}=\begin{bmatrix}
    \frac{bc}{a}+X & c & b\\
    c & \frac{ac}{b} & a\\
    b & a & \frac{ab}{c}
    \end{bmatrix}.
\end{equation}
If $\lambda=0$ is not the second eigenvalue of $\bar{\M}$, then it must be the largest eigenvalue. Since this would imply all eigenvalues of $\bar{\M}$ are non-positive, it follows that $\bar{\M}$ is negative semi-definite, or, equivalently, $-\bar \M$ is positive semi-definite. The latter implies that all of the leading principal minors of $-\bar \M$, i.e., the determinants of all upper-left $k$-by-$k$ sub-matrices of $-\bar{\M}$, are non-negative. 
In particular, the two-by-two leading principal minor of $-\bar{\M}$ gives
\begin{align*}
    \left(\frac{bc}{a}+X\right)\frac{ac}{b}-c^2\geq0,
\end{align*}
Since $ac/b<0$, this implies that
\begin{align}
    X\leq 0.
\end{align}
Conversely, when $X\leq0$, the remaining principal minors of $-\bar{\M}$ are also non-negative, so $-\bar{\M}$ is positive semi-definite. Recalling that $X\neq0$ on the axis under consideration, $\lambda=0$ therefore corresponds to the largest eigenvalue of $\bar{\M}$, with eigenvector $[0,-b,c]$, if and only if $X<0$.

\begin{figure}[htp]
    \centering
    \includegraphics[width=10cm]{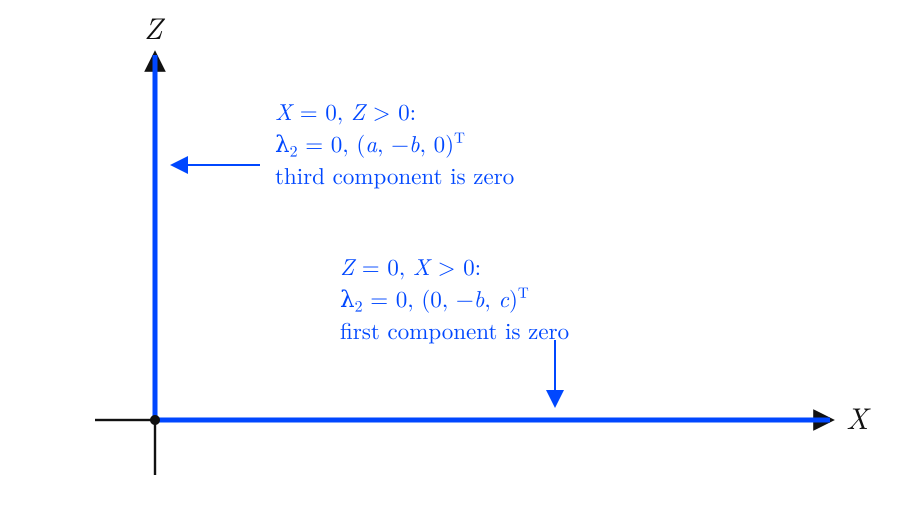}
    \caption{The second eigenvector, corresponding to a simple second eigenvalue $\lambda=0$, has a zero in the first position on the positive $X$-axis or a zero in the third position on the positive $Z$-axis as shown in blue.}\label{lambda_zeros}
\end{figure}

Consequently, $[0,-b,c]$ is the second eigenvector with eigenvalue $\lambda=0$ if and only if $X>0$ and $Z=0$. By repeating the above arguments, we find that on $\{Z\neq 0\}\cap\{X=0\}$, $[a,-b,0]$ is the second eigenvector if and only if $Z>0$. Pictorially, the second eigenvector corresponding to a simple second eigenvalue $\lambda=0$ has a zero component in the places indicated in Figure \ref{lambda_zeros}.\\

Now we must work out what happens to the zero component of the second eigenvector when we cross the lines $\{Z=0\}\cap\{X>0\}$ and $\{X=0\}\cap\{Z>0\}$ in each of the above two cases. As before, only one set of calculations is presented in detail since they are essentially the same. Using the same eigenvector $\bv=[0,-b,c]$ corresponding to the second eigenvalue $\lambda=0$ as previously,
we want to calculate the vector derivative of $\bv$ with respect to the scalar variable $Z$:
\begin{align}\label{v_derivative}
    \bv':=\frac{\partial \bv}{\partial Z}=  \begin{bmatrix}
                                            \frac{\partial x}{\partial Z}\\
                                            \frac{\partial y}{\partial Z}\\
                                            \frac{\partial z}{\partial Z}
                                        \end{bmatrix}
                                    =:   \begin{bmatrix}
                                            x'\\
                                            y'\\
                                            z'
                                        \end{bmatrix},
\end{align}
where $'$ denotes partial differentiation with respect to $Z$, and then extract the first component $x'$. To this end, we differentiate the eigenvalue equation (\ref{evalue_equation}) with respect $Z$, giving
\begin{align}\label{evalue_equation_derivative}
    0=\bar{\M}\bv'+(\bar{\M}'-\lambda ')\bv+\bv'\lambda=\bar{\M}\bv'+(\bar{\M}'-\lambda ')\bv,
\end{align}
where we have used the fact that $\lambda=0$ in the second equality. The matrix $\bar{\M}$, after replacing $U,V,W$ by equations (\ref{transformed_V}) and (\ref{transformed_U_W}), is
\begin{align}
    \bar{\M}=\begin{bmatrix}
                \frac{bc}{a}+X & c & b\\
                c & \frac{ac}{b} & a\\
                b & a & \frac{ab}{c}+Z
            \end{bmatrix},
\end{align}
so $\bar\M'=\frac{\partial \bar{\M}}{\partial Z}$ is given by
\begin{align}\label{matrix_derivative}
    \bar{\M}'=   \begin{bmatrix}
                    0 & 0 & 0\\
                    0 & 0 & 0\\
                    0 & 0 & 1\\
                \end{bmatrix}.
\end{align}
We can calculate the partial derivative $\lambda'=\frac{\partial\lambda}{\partial Z}$ of $\lambda$ with respect to $Z$ by pre-multiplying equation (\ref{evalue_equation_derivative}) by $\bv^T$ and rearranging the result for $\lambda'$ as follows:
\begin{align}\label{lambda_derivative}
0&=\bv^T\bar\M\bv'+\bv^T\bar\M'\bv-\lambda'\|\bv\|^2\nonumber\\
&=c^2-\lambda'(b^2+c^2),
\end{align}
so
\begin{align}\label{lambda_derivative2}
    \lambda'&=\frac{c^2}{b^2+c^2}=c^2, \qquad b^2+c^2=1.
\end{align}
The fact that $\bv$ is an eigenvector of $\bar{\M}$, with $\|\bv\|=1$, and $\bar{\M}$ is symmetric so $\bv^T\bar \M=\bv^T\bar \M^T=(\bar\M\bv)^T=(\lambda \bv)^T=(0.\bv)^T=0$ was applied in the first line of equation (\ref{lambda_derivative}). Since our interest is in $x'$, we can repeat the above calculation, but now pre-multiply equation (\ref{evalue_equation_derivative}) by $[1,0,0]$:
\begin{align}\label{pre_equation_x'}
    0&= \begin{bmatrix}
        1 & 0 & 0
        \end{bmatrix}\bar{\M}  \begin{bmatrix}
                    x'\\
                    y'\\
                    z'
                \end{bmatrix}+  \begin{bmatrix}
                                    1 & 0 & 0
                                \end{bmatrix}\bar{\M}'   \begin{bmatrix}
                                                            0\\
                                                            -b\\
                                                            c
                                                        \end{bmatrix}-  \begin{bmatrix}
                                                                            1 & 0 & 0
                                                                        \end{bmatrix}\lambda'   \begin{bmatrix}
                                                                                                    0\\
                                                                                                    -b\\
                                                                                                    c
                                                                                                \end{bmatrix}\nonumber\\
    &=  \begin{bmatrix}
            \Big(\frac{bc}{a}+X\Big) & 0 & 0
        \end{bmatrix}   \begin{bmatrix}
                            x'\\
                            y'\\
                            z'
                        \end{bmatrix}\nonumber\\
    &=  \begin{bmatrix}
            \frac{bc}{a} & c & b
        \end{bmatrix}   \begin{bmatrix}
                            x'\\
                            y'\\
                            z\
                        \end{bmatrix}+  \begin{bmatrix}
                                            X & 0 & 0
                                        \end{bmatrix}   \begin{bmatrix}
                                                            x'\\
                                                            y'\\
                                                            z'    
                                                        \end{bmatrix}\nonumber\\
    &=  \begin{bmatrix}
            \frac{bc}{a} & c & b
        \end{bmatrix}   \begin{bmatrix}
                            x'\\
                            y'\\
                            z\
                        \end{bmatrix}+Xx',
\end{align}
where the middle term  in the first line evaluates to zero. Now we just need an expression for the first term in the last line of the above computation (\ref{pre_equation_x'}). This can be achieved by choosing the vector we pre-multiply equation (\ref{evalue_equation_derivative}) by to be perpendicular to the eigenvector $[0,-b,c]$:
\begin{align*}
     0&=\begin{bmatrix}
            0 & c & b
        \end{bmatrix}\bar{\M}  \begin{bmatrix}
                    x'\\
                    y'\\
                    z'
                \end{bmatrix}+  \begin{bmatrix}
                                    0 & c & b
                                \end{bmatrix}\bar{\M}'   \begin{bmatrix}
                                                            0\\
                                                            -b\\
                                                            c
                                                        \end{bmatrix}-  \begin{bmatrix}
                                                                            0 & c & b
                                                                        \end{bmatrix}\lambda'   \begin{bmatrix}
                                                                                                    0\\
                                                                                                    -b\\
                                                                                                    c
                                                                                                \end{bmatrix}\nonumber\\
    &=  \begin{bmatrix}
            c^2+b^2 & \frac{ac^2}{b}+ab & ac+\frac{ab^2}{c}
        \end{bmatrix}   \begin{bmatrix}
                            x'\\
                            y'\\
                            z'
                        \end{bmatrix}+bc\nonumber\\
    &=\frac{a(c^2+b^2)}{bc}\left(\begin{bmatrix}
                                    \frac{bc}{a} & c & b
                                \end{bmatrix}   \begin{bmatrix}
                                                    x'\\
                                                    y'\\
                                                    z'
                                                \end{bmatrix}+\frac{b^2c^2}{a(b^2+c^2)}\right),
\end{align*}
so
\begin{align}\label{pre_equation_x'_2}
\begin{bmatrix}
    \frac{bc}{a} & c & b
\end{bmatrix}   \begin{bmatrix}
                    x'\\
                    y'\\
                    z'
                \end{bmatrix}&=-\frac{b^2c^2}{a(b^2+c^2)}.
\end{align}
Substituting equation (\ref{pre_equation_x'_2}) into equation (\ref{pre_equation_x'}) and rearranging for $x'$ gives
\begin{align}
    x'=\frac{\partial x}{\partial Z}=\frac{1}{X}\frac{b^2c^2}{a(b^2+c^2)}=\frac{1}{X}\frac{b^2c^2}{a}<0,
\end{align}
assuming $b^2+c^2=1$, $a,b,c<0$, and $X>0$. That is, the change in $x$ is proportional to and opposite in sign to the change in $Z$.

\begin{figure}[htp]
    \centering
    \includegraphics[width=14cm]{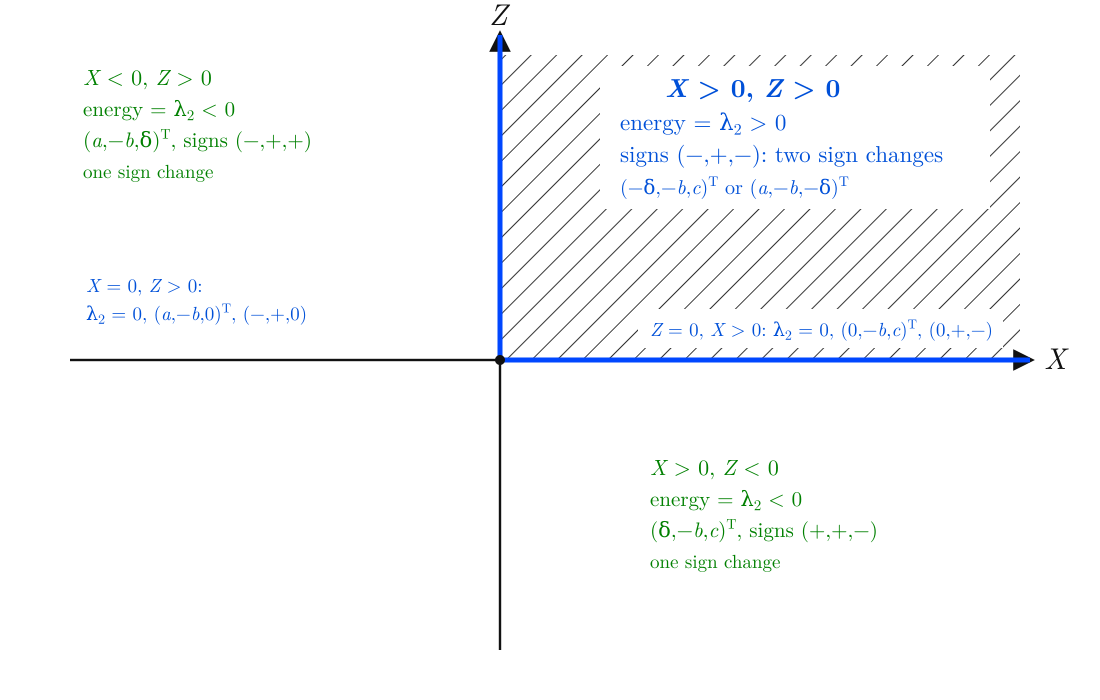}
    \caption{The second eigenvector has two sign changes and positive energy in the first quadrant $\{X>0\}\cap\{Z>0\}$.}\label{lambda_sign_change}
\end{figure}

Consequently, as we cross the line $\{Z=0\}\cap\{X>0\}$ into the $\{Z>0\}\cap\{X>0\}$ quadrant, the zero component of $(0,-b,c)$ decreases in value, so the resulting eigenvector has signs $(-,+,-)$ (since $b,c<0$ by assumption), which gives two sign changes. Conversely, when we cross $\{Z=0\}\cap\{X>0\}$ into the $\{Z<0\}\cap\{X>0\}$ quadrant, the zero component of $(0,-b,c)$ increases in value, so the resulting eigenvector has signs $(+,+,-)$, yielding one sign change. A similar analysis shows that if $(a,-b,0)$ is the second eigenvector on the set $\{X=0\}\cap\{Z>0\}$, the zero component has two sign changes in the $\{X>0\}\cap\{Z>0\}$ quadrant and one sign change in the $\{X<0\}\cap\{Z>0\}$ quadrant. See Figure \ref{lambda_sign_change}.\\ 

Finally, we check how the energies compare when the eigenvector has one or two sign changes in each of the above two cases. Recall that the energy is the corresponding eigenvalue when the vector in equation (\ref{energy_putative_counterexample}) is an eigenvector. Since the calculation is the same in both cases, we only provide details for $(0,-b,c)$. 

From equation (\ref{lambda_derivative2}), we know that $\lambda'=\frac{\partial\lambda}{\partial Z}$ is a positive constant. This means that as we increase (decrease) $Z$, $\lambda$ also increases (decreases). Since $x'<0$, the change in the zero component has the opposite sign to the change in $\lambda$. Thus the energy corresponding to $(\delta,-b,c)$, with $\delta>0$ and signs $(+,+,-)$, is negative; while the energy for $(-\delta,-b,c)$, with signs $(-,+,-)$, is positive. This shows that we have found a situation where the second eigenvector has two sign changes, namely, when $X>0$ and $Z>0$. In this regime the assumptions $|a|,|c|>|b|$ imply that the middle diagonal value satisfies $V<U,W$. The same conclusion holds when we vary the last component of $(a,-b,0)$ on $\{X=0\}\cap\{Z>0\}$.

Thus, $\bar{\M}$ has a second eigenvector with two sign changes, corresponding to a positive second eigenvalue that is simple, when $V<U,W$. Figure \ref{lambda_sign_change} summarises these results.\\ 

On the basis of the above approximation to the infinite potential well eigenvalue problem (\ref{simple_evalue_problem_perturbed_frac_Lap_1d}), we produce a counterexample to its second eigenfunction having only one zero in the next section.

\section{Infinite potential well counterexample}\label{sec_limiting_counterexample}
Fix $\epsilon$ small and positive, and assume that the subintervals $I_i=(x_i-\epsilon, x_i+\epsilon)$ are mutually disjoint in $(-1,1)$. Our aim in this section is to construct a counterexample to the second eigenfunction having only one sign change of eigenvalue problem (\ref{simple_evalue_problem_perturbed_frac_Lap_1d}), reproduced below for convenience:
\begin{equation}
 \left\{\begin{aligned}
        ((-\Delta)^s+\epsilon V_i)u(x)&=\lambda u(x), &x\in I_i, i=1,\ldots,k\nonumber\\
        u(x)&=0, &x\in \R\setminus \bigcup_{i=1}^k I_i\nonumber.
       \end{aligned}
 \right.
\end{equation}
As mentioned in the counterexample outline in Section \ref{ch_4_background}, the key to achieving this goal is to reduce problem (\ref{simple_evalue_problem_perturbed_frac_Lap_1d}) to the finite-dimensional eigenvalue problem solved in Section \ref{putative_counterexample}, by showing it is amenable to a regular perturbation analysis. We first review regular perturbation theory within the context of problem (\ref{simple_evalue_problem_perturbed_frac_Lap_1d}).

\subsection{Regular perturbation theory}\label{Kato_Rellich_theory}
Let $H_0$ be an operator with a discrete eigenvalue $\lambda_0$ that may be simple or degenerate. Consider the perturbed operator $T(\beta)=H_0+f(\beta)$, where the perturbation $f(\beta)$ is another operator that depends on a small perturbation parameter $\beta$ taking values in a connected region $R$ of the complex plane. The Kato-Rellich theory provides simple criteria under which $T(\beta)$ is guaranteed to have eigenvalues $\lambda(\beta)$ near $\lambda_0$ in the form of convergent power series in $\beta$ near zero:
\begin{align}\label{series_solution_reg_pert_th}
    \lambda(\beta)&=\lambda^0+\beta \lambda^1+\beta^2\lambda^2+\ldots\\
    \phi(\beta)&=\phi^0+\beta\phi^1+\beta^2\phi^2+\ldots,
\end{align}
where $\lambda(\beta)$ and $\phi(\beta)$ are the perturbed eigenvalue and corresponding perturbed eigenfunction, respectively. Here $\lambda^0=\lambda(0)=\lambda_0$ and
$\phi^0=\phi(0)$ is an eigenfunction of $H_0$ corresponding to
$\lambda_0$, while, for $n\geq1$, $\lambda^n$ and $\phi^n$ are the
$n$th-order coefficients, or corrections, in the respective expansions.

We are interested in the case where $\lambda_0$ is degenerate. In this case, there are two criteria to be satisfied: $T(\beta)$ needs to be an analytic family in the sense of Kato for $\beta$ near zero (see below); and $T(\beta)$ must be self-adjoint for each real-valued $\beta$ in $R$.
For the infinite potential well eigenvalue problem (\ref{simple_evalue_problem_perturbed_frac_Lap_1d}), the self-adjointness property follows naturally from the self-adjointness of the fractional Dirichlet Laplace operator on a bounded domain, as we show below. So the central difficulty to overcome is showing the operator $(-\Delta)^s+\epsilon V_i$, for $i=1,2,\ldots,k$, or a rescaled version of it (see below), is an analytic family in the sense of Kato. 
\begin{definition}[Analytic family of type (A)]\label{def_analytic_family_A} Let $R$ be a connected domain in $\mathbb{C}$, and let $T(\beta)$ be a closed operator with non-empty resolvent set for each $\beta$ in $R$. We say that $T(\beta)$ is an analytic family of type (A) if and only if the following two conditions are satisfied:
\begin{enumerate}
    \item[(i)] The operator domain of $T(\beta)$ is a set $D$ independent of $\beta$.
    \item[(ii)] For each $\phi\in D$, $T(\beta)\phi$ is a vector-valued analytic function of $\beta$.
\end{enumerate}
\end{definition}
Any analytic family of type (A) is also analytic in the sense of Kato. \cite[p.16, Vol. 4]{RS} So provided we can satisfy Definition \ref{def_analytic_family_A}, the Kato-Rellich theory of regular perturbations for degenerate eigenvalues guarantees the existence of perturbed eigenvalues (and corresponding eigenfunctions) analytic in the perturbation parameter. This is formalised in the theorem to follow, taken from \cite[Ch. XII.2, Vol. 4]{RS}. 
\begin{theorem}[Kato-Rellich theory for degenerate eigenvalues]\label{Kato_Rellich_degenerate_evalue} Assume the following:
\begin{itemize}
    \item[(i.)] $T(\beta)$ is an analytic family in the sense of Kato for $\beta$ near zero.
    \item[(ii).] $T(\beta)$ is self-adjoint for each $\beta$ real.
    \item[(iii.)] $\lambda_0$ is a discrete eigenvalue of multiplicity $k$ of the unperturbed operator $T(0)$. 
\end{itemize} 
Then there are $k$ not necessarily distinct single-valued functions, analytic near $\beta=0$, denoted by $\lambda^{(1)}(\beta),\lambda^{(2)}(\beta),\ldots,\lambda^{(k)}(\beta)$, with $\lambda^1(0)=...=\lambda^k(0)=\lambda_0$, so that $\lambda^{(i)}(\beta)$ 
are eigenvalues of $T(\beta)$ for $\beta$ near zero. Furthermore, these are the only eigenvalues near $\lambda_0$. 
\end{theorem} 

As we shall see in the final part of this section, Theorem \ref{Kato_Rellich_degenerate_evalue} reduces an infinite-dimensional degenerate eigenvalue problem to a finite-dimensional one, so that the existence of eigenvectors in the finite-dimensional case implies their existence in the infinite-dimensional problem. We now turn to applying the Kato-Rellich theory to construct a counterexample to the second eigenfunction having only one zero of eigenvalue problem (\ref{simple_evalue_problem_perturbed_frac_Lap_1d}).

\subsection{Application of the Kato-Rellich theory}\label{sec_KatoRellich_app}
For $\epsilon>0$, set $\Omega_\epsilon:=\bigcup_{i=1}^k I_i$, $I_i=(x_i-\epsilon,x_i+\epsilon)$, and let
\[
    T(\epsilon)
    :=(-\Delta)^s_{\mathrm{res},\Omega_\epsilon}
      +\epsilon M_V,
    \qquad
    D(T(\epsilon))
    =D\bigl((-\Delta)^s_{\mathrm{res},\Omega_\epsilon}\bigr),
\]
denote the operator associated with eigenvalue problem
\eqref{simple_evalue_problem_perturbed_frac_Lap_1d} on
$L^2(\Omega_\epsilon)$, where $M_V$ denotes multiplication by the bounded potential satisfying $V=V_i$ on $I_i$. This family $\{T(\epsilon)\}_{\epsilon>0}$ is not directly amenable to the degenerate regular perturbation analysis of Theorem \ref{Kato_Rellich_degenerate_evalue}.

First, the intervals $I_i$ varying with respect to $\epsilon$ implies so do the underlying space $L^2\left(\bigcup_{i=1}^k I_i\right)$
and the domain of the restricted fractional Laplacian. 
In its original form, the family is not an analytic family of operators on a fixed Hilbert space with a common domain.

Second, the original family does not provide a suitable unperturbed
operator with a ground-state eigenvalue of multiplicity $k$. As
$\epsilon$ tends to zero, the intervals collapse to the points
$x_1,\ldots,x_k$, so setting $\epsilon=0$ does not produce a
non-trivial limiting operator satisfying condition \textup{(iii.)} of
Theorem \ref{Kato_Rellich_degenerate_evalue}. On the other hand, for every fixed $\epsilon>0$, the restricted fractional Laplacian on $\Omega_\epsilon$
is not the direct sum of the restricted fractional Laplacians on the
individual intervals. Indeed, let $\phi_j$ be the positive
ground-state eigenfunction on $I_j$, extended by zero to $\R$. For
$x\in I_i$, where $i\neq j$, we have $\phi_j(x)=0$, and hence
\[
\begin{aligned}
    (-\Delta)^s_{\mathrm{res},\Omega_\epsilon}\phi_j(x)
    &=
    -c_{1,s}\int_{I_j}
    \frac{\phi_j(y)}{|x-y|^{1+2s}}\,dy\\
    &<0
    =
    \lambda_1(I_j)\phi_j(x).
\end{aligned}
\]
Thus the zero extension of $\phi_j$ is not an eigenfunction of the
operator on $\Omega_\epsilon$, as the disjoint intervals remain coupled
through the non-local interaction terms. In particular, the $k$
one-interval ground states do not give a $k$-dimensional eigenspace on
$\Omega_\epsilon$. By the standard simplicity result for the
first eigenvalue of the restricted fractional Laplacian, the actual
ground-state eigenvalue on $\Omega_\epsilon$ is simple. The necessary
$k$-fold degeneracy is thus absent from the original family.

Finally, the spectral scale becomes singular as $\epsilon$ tends to
zero. Let $I:=(-1,1)$.
The restricted fractional Laplacian on each interval $I_i$ has
discrete eigenvalues satisfying
\[
    \lambda_m(I_i)
    =
    \epsilon^{-2s}\lambda_m(I),
\]
so for each fixed $m$, $\lambda_m(I_i)$ diverges at the rate
$\epsilon^{-2s}$ as $\epsilon\to0^+$.

These considerations suggest rescaling the original family of operators in a way that preserves its $L^2$ Hilbert space structure,
while removing the singular eigenvalue scale. This allows us to obtain an equivalent spectral problem on a fixed Hilbert space that is amenable to degenerate regular perturbation analysis.\\

We begin by identifying each interval $I_j$ with the fixed reference
interval $I=(-1,1)$ through the transformation
\begin{align}\label{rescale_pert_frac_Lap_Dir}
(x_j-\epsilon,x_j+\epsilon)\ni y
\longmapsto
\frac{y-x_j}{\epsilon}\in I,
\end{align}
where $\epsilon>0$ and $j=1,\ldots,k$. At the level of function
spaces, the corresponding change of variables is used to define a unitary
identification of the varying spaces $L^2(\Omega_\epsilon)$ with the
fixed Hilbert space
\[
\mathcal H:=\prod_{j=1}^kL^2(I),
\]
equipped with the inner product
\begin{align}\label{prod_sp_IP}
\langle \fvect{u},\fvect{v}\rangle_{\mathcal H}
:=
\sum_{j=1}^k
\langle u_j,v_j\rangle_{L^2(I)}.
\end{align}
Since the product is finite and each $L^2(I)$ is complete,
$\mathcal H$ is a Hilbert space. 

Now determine this unitary map. Fix $\epsilon>0$, and define the linear map
\begin{align}\label{unitary_rescaling_map}
    \mathcal U_\epsilon:L^2(\Omega_\epsilon)
    &\longrightarrow\mathcal H,\nonumber\\
    [\mathcal U_\epsilon u]_i(z)
    &=:u_i(z)
      =\sqrt{\epsilon}\,u(x_i+\epsilon z),
      \qquad z\in I.
\end{align}
For any
$u,v\in L^2(\Omega_\epsilon)$, the change of variables
$x=x_i+\epsilon z$, $dx=\epsilon\,dz$, gives
\begin{align*}
    \langle\mathcal U_\epsilon u,\mathcal U_\epsilon v
    \rangle_{\mathcal H}
    &=
    \sum_{i=1}^k
    \int_I
    \epsilon\,
    u(x_i+\epsilon z)
    v(x_i+\epsilon z)\,dz\\
    &=
    \sum_{i=1}^k
    \int_{I_i}u(x)v(x)\,dx\\
    &=
    \langle u,v\rangle_{L^2(\Omega_\epsilon)}.
\end{align*}
Thus $\mathcal U_\epsilon$ preserves the inner product and, in
particular, the norm. To show that $\mathcal U_\epsilon$ is onto, let
$\fvect{f}=(f_1,\ldots,f_k)\in\mathcal H$ and define 
$\mathcal V_\epsilon\fvect{f}$
on $\Omega_\epsilon$
by
\begin{align}\label{inverse_unitary_rescaling_map}
    [\mathcal V_\epsilon\fvect{f}](x)
    :=
    \frac{1}{\sqrt{\epsilon}}
    f_i\left(\frac{x-x_i}{\epsilon}\right),
    \qquad x\in I_i,\quad i=1,\ldots,k.
\end{align}
As the intervals $I_i$ are pairwise disjoint, every
$x\in\Omega_\epsilon$ belongs to a unique interval $I_i$. Moreover,
$x\in I_i$ implies that $\frac{x-x_i}{\epsilon}\in I$.
Hence \eqref{inverse_unitary_rescaling_map} defines a piecewise
function on $\Omega_\epsilon$. Furthermore,
\begin{align*}
    \|\mathcal V_\epsilon\fvect{f}\|_{L^2(\Omega_\epsilon)}^2
    &=
    \sum_{i=1}^k
    \int_{I_i}
    \frac{1}{\epsilon}
    \left|
    f_i\left(\frac{x-x_i}{\epsilon}\right)
    \right|^2\,dx\\
    &=
    \sum_{i=1}^k\int_I|f_i(z)|^2\,dz\\
    &=
    \|\fvect{f}\|_{\mathcal H}^2,
\end{align*}
so $\mathcal V_\epsilon\fvect{f}\in L^2(\Omega_\epsilon)$.
Direct substitution into (\ref{unitary_rescaling_map}) gives
\begin{align*}
    [\mathcal U_\epsilon\mathcal V_\epsilon\fvect{f}]_i(z)
    &=
    \sqrt{\epsilon}\,
    [\mathcal V_\epsilon\fvect{f}](x_i+\epsilon z)\\
    &=
    \sqrt{\epsilon}\,
    \frac{1}{\sqrt{\epsilon}}
    f_i\left(
        \frac{x_i+\epsilon z-x_i}{\epsilon}
    \right)\\
    &=
    f_i(z).
\end{align*}
Therefore,
\[
    \mathcal U_\epsilon\mathcal V_\epsilon\fvect{f}
    =
    \fvect{f}.
\]
So $\mathcal U_\epsilon$ is onto.
As it preserves the norm, it is also
one-to-one. 
It follows that $\mathcal U_\epsilon$ is unitary, and
\[
    \mathcal V_\epsilon=\mathcal U_\epsilon^{-1}.
\]
The factor $\sqrt{\epsilon}$ in
\eqref{unitary_rescaling_map} is chosen so that its square cancels the
factor $1/\epsilon$ arising from the change of variables
$dz=dx/\epsilon$, thereby preserving the norm.

For $\epsilon>0$, the unitary map $\mathcal U_\epsilon$ transports
$T(\epsilon)$ from the varying space $L^2(\Omega_\epsilon)$ to the
fixed Hilbert space $\mathcal H$. Thus the original eigenvalue problem
\[
    T(\epsilon)u=\lambda_\epsilon u
\]
is equivalent, under $\fvect{u}=\mathcal U_\epsilon u$,
to
\[
    \mathcal U_\epsilon T(\epsilon)\mathcal U_\epsilon^{-1}
    \fvect{u}
    =
    \lambda_\epsilon\fvect{u}.
\]

The unitary transformation alone, however, does not normalise the singular spectral
scale associated with the shrinking intervals. Recall that the
eigenvalues of the one-interval restricted fractional Laplacian
satisfy
\begin{align}\label{singular_scale_factor}
    \lambda_m(I_i)
    =
    \epsilon^{-2s}\lambda_m(I),
\end{align}
which suggests multiplying the transformed
eigenvalue problem by $\epsilon^{2s}$. Motivated by this one-interval scaling, we define
\begin{equation}\label{definition_rescaled_operator}
    \overline T(\epsilon)
    :=
    \epsilon^{2s}
    \mathcal U_\epsilon
    T(\epsilon)
    \mathcal U_\epsilon^{-1},
    \qquad \epsilon>0.
\end{equation}
As we show below, this normalisation also produces
a finite, nontrivial limiting fractional Laplace operator for the full
nonlocal problem on $\Omega_\epsilon$.

We now derive the potential and restricted fractional Laplace terms
in the re-scaled operator $\overline T(\epsilon)$. First, let
$\fvect{u}=(u_1,\ldots,u_k)\in\mathcal H$ and set
\[
    u:=\mathcal U_\epsilon^{-1}\fvect{u}.
\]
For $z\in I$, we have $x_i+\epsilon z\in I_i$, and hence
\begin{align*}
    [\mathcal U_\epsilon M_V
      \mathcal U_\epsilon^{-1}\fvect{u}]_i(z)
    &=
    \sqrt{\epsilon}\,
    [M_Vu](x_i+\epsilon z)\\
    &=
    \sqrt{\epsilon}\,
    V(x_i+\epsilon z)u(x_i+\epsilon z)\\
    &=
    \sqrt{\epsilon}\,
    V_i
    \frac{1}{\sqrt{\epsilon}}
    u_i\left(
        \frac{x_i+\epsilon z-x_i}{\epsilon}
    \right)\\
    &=
    V_i u_i(z).
\end{align*}
Therefore,
\[
    \mathcal U_\epsilon M_V\mathcal U_\epsilon^{-1}
    =
    \overline V,
\]
where
\[
    \overline V\fvect{u}
    :=
    (V_1u_1,\ldots,V_ku_k).
\]
Since the potential term in the original operator is $\epsilon M_V$,
its contribution to the normalised re-scaled operator is
\[
    \epsilon^{2s}
    \mathcal U_\epsilon(\epsilon M_V)
    \mathcal U_\epsilon^{-1}
    =
    \epsilon^{1+2s}\overline V.
\]

Next, we re-scale the restricted fractional Laplace term. 
Fix $x=x_i+\epsilon z\in I_i$. The restricted fractional Laplace operator on $\Omega_\epsilon$
may be decomposed as
\begin{align}
(-\Delta)^s_{\mathrm{res},\Omega_\epsilon}u(x)
&=
c_s\pv\int_{I_i}
\frac{u(x)-u(y)}{|x-y|^{1+2s}}\,dy\nonumber\\
&\quad+
c_s\sum_{j\neq i}\int_{I_j}
\frac{u(x)-u(y)}{|x-y|^{1+2s}}\,dy\nonumber\\
&\quad+
c_s\int_{\R\setminus\Omega_\epsilon}
\frac{u(x)}{|x-y|^{1+2s}}\,dy.
\label{full_union_operator_split}
\end{align}
Let $u^{(i)}$ denote the restriction of $u$ to $I_i$, extended by zero
on $\R\setminus I_i$. The corresponding one-interval restricted
operator satisfies
\begin{align}
(-\Delta)^s_{\mathrm{res},I_i}u^{(i)}(x)
&=
c_s\pv\int_{I_i}
\frac{u(x)-u(y)}{|x-y|^{1+2s}}\,dy
+
c_s\int_{\R\setminus I_i}
\frac{u(x)}{|x-y|^{1+2s}}\,dy\nonumber\\
&=
c_s\pv\int_{I_i}
\frac{u(x)-u(y)}{|x-y|^{1+2s}}\,dy\nonumber\\
&\quad+
c_s\sum_{j\neq i}\int_{I_j}
\frac{u(x)}{|x-y|^{1+2s}}\,dy
+
c_s\int_{\R\setminus\Omega_\epsilon}
\frac{u(x)}{|x-y|^{1+2s}}\,dy,
\label{single_interval_operator_split}
\end{align}
where
\[
    \R\setminus I_i
    =
    (\R\setminus\Omega_\epsilon)
    \cup\bigcup_{j\neq i}I_j.
\]
Comparing \eqref{full_union_operator_split} and
\eqref{single_interval_operator_split} leads to the identity
\begin{align}
(-\Delta)^s_{\mathrm{res},\Omega_\epsilon}u(x)
&=
(-\Delta)^s_{\mathrm{res},I_i}u^{(i)}(x)
-
c_s\sum_{j\neq i}\int_{I_j}
\frac{u(y)}{|x-y|^{1+2s}}\,dy.
\label{union_operator_single_interval_identity}
\end{align}
This identity highlights the nonlocal interaction between
distinct intervals.

Recall from \eqref{unitary_rescaling_map} that
\[
    u(x_i+\epsilon z)
    =
    \frac{u_i(z)}{\sqrt{\epsilon}}.
\]
Extend $u_i$ by zero on $\R\setminus I$. In the first term of
\eqref{union_operator_single_interval_identity}, make the change of
variables $y=x_i+\epsilon w$, so that $dy=\epsilon\,dw$. Then
\begin{align}
(-\Delta)^s_{\mathrm{res},I_i}u^{(i)}(x_i+\epsilon z)
&=
c_s\pv\int_{\R}
\frac{\epsilon^{-1/2}(u_i(z)-u_i(w))}
{\epsilon^{1+2s}|z-w|^{1+2s}}\,
\epsilon\,dw\nonumber\\
&=
\epsilon^{-1/2-2s}
(-\Delta)^s_{\mathrm{res},I}u_i(z).
\label{rescaled_single_interval_term}
\end{align}
Consequently,
\begin{align}
\epsilon^{2s}\sqrt{\epsilon}\,
(-\Delta)^s_{\mathrm{res},I_i}
u^{(i)}(x_i+\epsilon z)
=
(-\Delta)^s_{\mathrm{res},I}u_i(z).
\label{rescaled_single_interval_term_multiplied}
\end{align}
For $j\neq i$, set $y=x_j+\epsilon w$. Since
\[
    u(x_j+\epsilon w)
    =
    \epsilon^{-1/2}u_j(w),
\]
we have
\begin{align}
\int_{I_j}
\frac{u(y)}
{|x_i+\epsilon z-y|^{1+2s}}\,dy
&=
\epsilon^{1/2}
\int_I
\frac{u_j(w)}
{|x_i-x_j+\epsilon(z-w)|^{1+2s}}\,dw.
\label{rescaled_interaction_integral}
\end{align}
Multiplying \eqref{union_operator_single_interval_identity} by
$\epsilon^{2s}\sqrt{\epsilon}$ and using
\eqref{rescaled_single_interval_term_multiplied} and
\eqref{rescaled_interaction_integral}, we obtain
\begin{align}\label{correct_rescaled_frac_lap}
&\epsilon^{2s}\sqrt{\epsilon}\,
(-\Delta)^s_{\mathrm{res},\Omega_\epsilon}
u(x_i+\epsilon z)\nonumber\\
&\qquad=
(-\Delta)^s_{\mathrm{res},I}u_i(z)
-
c_s\epsilon^{1+2s}
\sum_{j\neq i}
\int_I
\frac{u_j(w)}
{|x_i-x_j+\epsilon(z-w)|^{1+2s}}\,dw.
\end{align}

Let $\lambda_\epsilon$ be an eigenvalue of the original operator and
suppose that
\[
    T(\epsilon)u=\lambda_\epsilon u.
\]
Set
\[
    \fvect{u}:=\mathcal U_\epsilon u,
\]
so that, by \eqref{unitary_rescaling_map},
\[
    u_i(z)
    =
    \sqrt{\epsilon}\,u(x_i+\epsilon z).
\]
Since
\[
    T(\epsilon)
    =
    (-\Delta)^s_{\mathrm{res},\Omega_\epsilon}
    +
    \epsilon M_V
\]
and $V=V_i$ on $I_i$, evaluating the original eigenvalue equation at
$x=x_i+\epsilon z$ gives
\[
    (-\Delta)^s_{\mathrm{res},\Omega_\epsilon}
    u(x_i+\epsilon z)
    +
    \epsilon V_i u(x_i+\epsilon z)
    =
    \lambda_\epsilon u(x_i+\epsilon z).
\]
Multiplying this equation by $\epsilon^{2s}\sqrt{\epsilon}$ yields
\begin{align*}
&\epsilon^{2s}\sqrt{\epsilon}\,
(-\Delta)^s_{\mathrm{res},\Omega_\epsilon}
u(x_i+\epsilon z)\\
&\qquad
+
\epsilon^{1+2s}V_i
\bigl[\sqrt{\epsilon}\,u(x_i+\epsilon z)\bigr]
=
\epsilon^{2s}\lambda_\epsilon
\bigl[\sqrt{\epsilon}\,u(x_i+\epsilon z)\bigr].
\end{align*}
Using
\[
    \sqrt{\epsilon}\,u(x_i+\epsilon z)=u_i(z)
\]
and the re-scaling identity \eqref{correct_rescaled_frac_lap}, we
obtain
\begin{align}\label{rescaled_evalue_problem}
    \mu_\epsilon u_i(z)
    &=
    (-\Delta)^s_{\mathrm{res},I}u_i(z)
    +
    \epsilon^{1+2s}V_i u_i(z)\nonumber\\
    &\quad-
    c_s\epsilon^{1+2s}
    \sum_{j\neq i}
    \int_I
    \frac{u_j(w)}
    {|x_i-x_j+\epsilon(z-w)|^{1+2s}}\,dw,
\end{align}
for $z\in I$ and $i=1,\ldots,k$, where
\[
    \mu_\epsilon
    :=
    \epsilon^{2s}\lambda_\epsilon.
\]
Thus \eqref{rescaled_evalue_problem} is the componentwise form of the
re-scaled eigenvalue equation
\[
    \overline T(\epsilon)\fvect{u}
    =
    \mu_\epsilon\fvect{u}.
\]
Indeed, using
\[
    \overline T(\epsilon)
    =
    \epsilon^{2s}
    \mathcal U_\epsilon
    T(\epsilon)
    \mathcal U_\epsilon^{-1},
\]
we have
\begin{align*}
    \overline T(\epsilon)\fvect{u}
    &=
    \epsilon^{2s}
    \mathcal U_\epsilon
    T(\epsilon)
    \mathcal U_\epsilon^{-1}\fvect{u}\\
    &=
    \epsilon^{2s}
    \mathcal U_\epsilon
    T(\epsilon)u\\
    &=
    \epsilon^{2s}\lambda_\epsilon
    \mathcal U_\epsilon u\\
    &=
    \mu_\epsilon\fvect{u}.
\end{align*}
We see the unitary map $\mathcal U_\epsilon$ relates the original and
re-scaled eigenfunctions, while the corresponding eigenvalues are
related by
\[
    \mu_\epsilon=\epsilon^{2s}\lambda_\epsilon.
\]
The calculation confirms that
the normalisation suggested by the one-interval scaling
\eqref{singular_scale_factor} removes the singular fractional
Laplace scale of the full nonlocal problem.

Let
\[
    A_I:=(-\Delta)^s_{\mathrm{res},I},
    \qquad
    D_I:=D(A_I),
\]
and define
\[
    D:=\prod_{i=1}^kD_I.
\]
The action of the re-scaled family is therefore
\begin{align}\label{rescaled_family_pert_op}
    [\overline T(\epsilon)\fvect{u}]_i(z)
    &=
    (-\Delta)^s_{\mathrm{res},I}u_i(z)
    +
    \epsilon^{1+2s}V_i u_i(z)\nonumber\\
    &\quad-
    c_s\epsilon^{1+2s}
    \sum_{j\neq i}
    \int_I
    \frac{u_j(w)}
    {|x_i-x_j+\epsilon(z-w)|^{1+2s}}\,dw.
\end{align}

As $d:=\min_{i\neq j}|x_i-x_j|>0$, and $z,w\in I=(-1,1)$, implies $|z-w|\leq 2$, we have
\[
    |x_i-x_j+\epsilon(z-w)|
    \geq |x_i-x_j|-\epsilon|z-w|
    \geq d-2\epsilon.
\]
It follows that for all sufficiently small $\epsilon>0$, there is a constant
$c>0$, independent of $i,j,z,w$, such that
\[
    |x_i-x_j+\epsilon(z-w)|\geq c
\]
whenever $i\neq j$ and $z,w\in I$. Consequently,
\[
    \frac{1}
    {|x_i-x_j+\epsilon(z-w)|^{1+2s}}
    \leq c^{-1-2s}.
\]
For $u_j\in L^2(I)$, Hölder's inequality gives
\begin{align*}
    \left\|
    \int_I
    \frac{u_j(w)}
    {|x_i-x_j+\epsilon(\,\cdot-w)|^{1+2s}}\,dw
    \right\|_{L^2(I)}
    &\leq
    c^{-1-2s}|I|\,
    \|u_j\|_{L^2(I)}.
\end{align*}
Since there are only finitely many pairs $i\neq j$, the full
interaction term defines a bounded operator on $\mathcal H$.

The potential operator
\[
    \overline V\fvect{u}
    =
    (V_1u_1,\ldots,V_ku_k)
\]
is also bounded on $\mathcal H$. Hence
\[
    \overline T(\epsilon)
    =
    \bigoplus_{i=1}^k(-\Delta)^s_{\mathrm{res},I}+B_\epsilon,
\]
where $B_\epsilon$ is bounded on $\mathcal H$. Since adding an
everywhere-defined bounded operator does not change the domain of an
unbounded operator, it follows that
\[
    D(\overline T(\epsilon))
    =
    \prod_{i=1}^kD(A_I)
    =
    D
\]
for all sufficiently small $\epsilon>0$.\\

In the re-scaled eigenvalue problem, $\epsilon>0$ is a real
parameter. The analyticity required in Theorem
\ref{Kato_Rellich_degenerate_evalue}, however, is analyticity of the
operator family on an open complex neighbourhood of the unperturbed
parameter value $0$. It is therefore not sufficient to consider the
dependence of the re-scaled operator only for positive real values of
$\epsilon$.

For $s=p/q\in(0,1)$ rational, we introduce the new parameter
\[
    \beta=\epsilon^{1/q},
    \qquad\text{so that}\qquad
    \epsilon=\beta^q.
\]
Then the perturbing terms in \eqref{rescaled_family_pert_op} containing
the factor
\[
    \epsilon^{1+2s}
    =
    \epsilon^{1+2p/q}
    =
    \beta^{q+2p}
\]
depend on the new parameter through integer powers of $\beta$. The
remaining dependence of the interaction kernel on $\epsilon$ is through
\[
    x_i-x_j+\epsilon(z-w)
    =
    x_i-x_j+\beta^q(z-w).
\]

Set
\[
    d_{ij}:=x_i-x_j,
    \qquad i\neq j.
\]
Since the centres $x_i$ are distinct and finite in number,
\[
    d_*:=\min_{i\neq j}|x_i-x_j|>0.
\]
As $z,w\in I=(-1,1)$, we also have $|z-w|<2$. Hence,
\[
    \left|
    \frac{\beta^q(z-w)}{d_{ij}}
    \right|
    <
    \frac{2|\beta|^q}{d_*}.
\]
Choose $\rho>0$ sufficiently small so that
\[
    2\rho^q<d_*,
\]
and set
\[
    \theta:=\frac{2\rho^q}{d_*}<1.
\]
Then
\[
    \left|
    \frac{\beta^q(z-w)}{d_{ij}}
    \right|
    <\theta<1
\]
for all $|\beta|<\rho$, $z,w\in I$, and $i\neq j$. This is
the condition under which the generalised binomial expansion used below
is valid.

To relate this expansion to the interaction kernel of the original
re-scaled problem, we first restrict to the real range
$0<\beta<\rho$. In this case,
\[
    |\beta^q(z-w)|<|d_{ij}|,
\]
so $d_{ij}+\beta^q(z-w)$ has the same sign as $d_{ij}$. It follows that
\[
    |d_{ij}+\beta^q(z-w)|
    =
    |d_{ij}|
    \left(
        1+\frac{\beta^q(z-w)}{d_{ij}}
    \right),
\]
where the factor in parentheses is positive. The generalised binomial theorem yields
\begin{align}\label{analytic_kernel_expansion}
\frac{1}{|d_{ij}+\beta^q(z-w)|^{1+2s}}
&=
\frac{1}{|d_{ij}|^{1+2s}}
\left(
1+\frac{\beta^q(z-w)}{d_{ij}}
\right)^{-(1+2s)}
\nonumber\\
&=
\frac{1}{|d_{ij}|^{1+2s}}
\sum_{n=0}^{\infty}
\binom{-(1+2s)}{n}
\left(
\frac{\beta^q(z-w)}{d_{ij}}
\right)^n.
\end{align}
Regarded as a power series in the complex
variable $\beta$, the right-hand side of
\eqref{analytic_kernel_expansion} defines an analytic function on the disc
$|\beta|<\rho$. For $0<\beta<\rho$, this function agrees with the
original interaction kernel and therefore defines its analytic extension
to the complex disc $|\beta|<\rho$:
\begin{align}\label{definition_analytic_kernel}
    K_{ij}(\beta;z,w)
    :=
    \frac{1}{|d_{ij}|^{1+2s}}
    \sum_{n=0}^{\infty}
    \binom{-(1+2s)}{n}
    \left(
        \frac{\beta^q(z-w)}{d_{ij}}
    \right)^n.
\end{align}

The operator family $\overline{T}(\beta)$, for complex
$|\beta|<\rho$, is finally given by
\begin{align}\label{rescaled_family_pert_op_Taylor}
    [\overline{T}(\beta)\fvect{u}]_i(z)
    &=
    (-\Delta)^{p/q}_{\mathrm{res},I}u_i(z)
    +\beta^{q+2p}V_i u_i(z)
    \nonumber\\
    &\quad
    -c_s\beta^{q+2p}
    \sum_{j\neq i}
    \int_I
    K_{ij}(\beta;z,w)u_j(w)\,dw.
\end{align}
When $\beta$ is real, positive, and sufficiently close to zero,
\eqref{rescaled_family_pert_op_Taylor} agrees with the re-scaled
operator \eqref{rescaled_family_pert_op} under the substitution
$\epsilon=\beta^q$. This family consists of closed operators with non-empty resolvent set on a sufficiently small connected neighbourhood of zero contained in $|\beta|<\rho$. It is an analytic family of type~(A) and self-adjoint for real $\beta$ on this neighbourhood. Moreover, the unperturbed operator has a discrete ground-state eigenvalue of multiplicity $k$. These properties are established in the three lemmas that follow.\\

\begin{lemma}
\label{rescaled_pert_op_nonempty_resolvent}
Let $s=p/q\in(0,1)$ be rational. There exists
$\rho_0\in(0,\rho]$ such that for every $\beta\in\mathbb{C}$ with
$|\beta|<\rho_0$, the operator $\overline{T}(\beta)$ given by
\eqref{rescaled_family_pert_op_Taylor} is closed on the fixed domain $D:=\prod_{i=1}^k D\bigl((-\Delta)^s_{\mathrm{res},I}\bigr)$
and has non-empty resolvent set. In fact,
\[
    i\in\rho\bigl(\overline{T}(\beta)\bigr)
\]
for all such $\beta$.
\end{lemma}

\begin{proof}
At $\beta=0$,
\[
    \overline{T}(0)
    =
    \bigoplus_{i=1}^k (-\Delta)^s_{\mathrm{res},I}
\]
on $D$. $(-\Delta)^s_{\mathrm{res},I}$ is self-adjoint on
$D((-\Delta)^s_{\mathrm{res},I})$ by
Lemma~\ref{lemma_self_adj_res_frac_Lapl}, so the finite direct sum
$\overline{T}(0)$ is self-adjoint on $D$. It is therefore also closed by the basic criterion for self-adjointness theorem.
\cite[Thm.~VIII.3, p.~257, Vol.~1]{RS} The same theorem implies that
$\overline{T}(0)-i$ is bijective. We now show that its inverse is bounded.

Let $u\in D$. We have
\[
    \|(\overline{T}(0)-iI)u\|_{\mathcal H}^2
    =
    \|\overline{T}(0)u\|_{\mathcal H}^2
    +
    \|u\|_{\mathcal H}^2
    \geq
    \|u\|_{\mathcal H}^2,
\]
and hence
\begin{equation}\label{self_adjoint_inequality}
    \|(\overline{T}(0)-iI)u\|_{\mathcal H}
    \geq
    \|u\|_{\mathcal H}.
\end{equation}
For any $f\in\mathcal H$, bijectivity gives a unique $u\in D$ such that
\[
    (\overline{T}(0)-iI)u=f,
    \qquad
    u=(\overline{T}(0)-iI)^{-1}f.
\]
Applying \eqref{self_adjoint_inequality} gives
\[
    \|(\overline{T}(0)-iI)^{-1}f\|_{\mathcal H}
    \leq
    \|f\|_{\mathcal H}.
\]
Thus
\[
    \|(\overline{T}(0)-i)^{-1}\|_{op}
    \leq 1,
\]
so the inverse is bounded. 

It follows by definition
\cite[p.~253, Vol.~1]{RS} that
\[
    i\in\rho(\overline{T}(0)).
\]

For $|\beta|<\rho$, write
\[
    \overline{T}(\beta)
    =
    \overline{T}(0)+B(\beta),
\]
where
\begin{align*}
    [B(\beta)\fvect{u}]_i(z)
    &=
    \beta^{q+2p}V_i u_i(z) \\
    &\quad
    -c_s\beta^{q+2p}
    \sum_{j\neq i}
    \int_I
    K_{ij}(\beta;z,w)u_j(w)\,dw .
\end{align*}
Using the estimates established earlier when deriving the analytic extension of the rescaled family of operators, the $n$th term in the series
\eqref{definition_analytic_kernel} defining $K_{ij}$ is bounded in
absolute value by
\[
    M_n
    :=
    \frac{1}{d_*^{1+2s}}
    \left|
        \binom{-(1+2s)}{n}
    \right|
    \theta^n
\]
uniformly for $|\beta|<\rho$, $z,w\in I$, and $i\neq j$. The generalised binomial theorem yields
\[
    \sum_{n=0}^{\infty} M_n
    =
    \frac{1}{d_*^{1+2s}}
    (1-\theta)^{-(1+2s)}
    <\infty,
\]
so the Weierstrass $M$-test implies that the series defining
$K_{ij}(\beta;z,w)$ converges absolutely and uniformly on this set. In particular,
\[
    |K_{ij}(\beta;z,w)|
    \leq
    \frac{1}{d_*^{1+2s}}
    (1-\theta)^{-(1+2s)}
\]
uniformly for $|\beta|<\rho$, $z,w\in I$, and $i\neq j$.

For $i\neq j$, define
\[
    [C_{ij}(\beta)u](z)
    :=
    \int_I K_{ij}(\beta;z,w)u(w)\,dw .
\]
If $M>0$ is a uniform bound for the kernels, then H\"older's inequality leads to
\[
    |C_{ij}(\beta)u(z)|
    \leq
    M\int_I |u(w)|\,dw
    \leq
    M|I|^{1/2}\|u\|_{L^2(I)}
\]
uniformly for $z\in I$. Therefore,
\[
    \|C_{ij}(\beta)u\|_{L^2(I)}
    \leq
    M|I|\|u\|_{L^2(I)},
\]
which shows that each $C_{ij}(\beta)$ is a bounded linear operator on $L^2(I)$. This bound is uniform for $|\beta|<\rho$.

As $V_i\in L^\infty(I)$ and the interaction operators
$C_{ij}(\beta)$ are uniformly bounded in operator norm for
$|\beta|<\rho$, each component satisfies
\[
\begin{aligned}
    \|[B(\beta)\mathbf u]_i\|_{L^2(I)}
    &\leq
    |\beta|^{q+2p}
    \left(
        \|V_i\|_{L^\infty(I)}\|u_i\|_{L^2(I)}
        +
        c_s\sum_{j\neq i}
        \|C_{ij}(\beta)\|_{\mathrm{op}}
        \|u_j\|_{L^2(I)}
    \right).
\end{aligned}
\]
There are finitely many components and interaction terms, so there
exists a constant $C>0$, independent of $\beta$, such that
\[
    \|B(\beta)\mathbf u\|_{\mathcal H}
    \leq
    C|\beta|^{q+2p}\|\mathbf u\|_{\mathcal H}
\]
for $|\beta|<\rho$. Thus $B(\beta)$ is a bounded operator on
$\mathcal H$ for every $|\beta|<\rho$.

Since $\overline T(0)$ is closed and $B(\beta)$ is bounded, we may
verify directly that $\overline T(\beta)=\overline T(0)+B(\beta)$ is
closed on $D$. Let $\mathbf u_n\in D$ satisfy
\[
    \mathbf u_n\to\mathbf u
    \qquad\text{and}\qquad
    \overline T(\beta)\mathbf u_n\to\mathbf f
\]
in $\mathcal H$. Since $B(\beta)$ is bounded, it is continuous, so
\[
    B(\beta)\mathbf u_n\to B(\beta)\mathbf u.
\]
Hence
\[
    \overline T(0)\mathbf u_n
    =
    \overline T(\beta)\mathbf u_n-B(\beta)\mathbf u_n
    \to
    \mathbf f-B(\beta)\mathbf u.
\]
Now $\overline T(0)$ is closed, so it follows that
\[
    \mathbf u\in D
    \qquad\text{and}\qquad
    \overline T(0)\mathbf u
    =
    \mathbf f-B(\beta)\mathbf u.
\]
Therefore,
\[
    \overline T(\beta)\mathbf u=\mathbf f,
\]
and $\overline T(\beta)$ is closed on $D$.

From earlier in this proof, we showed that
\[
    \|B(\beta)\|_{\mathrm{op}}
    \leq
    C|\beta|^{q+2p},
\]
implying 
\[
    \|B(\beta)\|_{\mathrm{op}}
    \longrightarrow 0
    \qquad\text{as }\beta\longrightarrow 0.
\]
This smallness will now be used to prove that $i$ remains in the
resolvent set of $\overline{T}(\beta)$ for all sufficiently small
$\beta$. 

We concluded previously that $i\in\rho\bigl(\overline{T}(0)\bigr)$
and the inverse $(\overline{T}(0)-iI)^{-1}$ is bounded:
\[
    \|(\overline{T}(0)-iI)^{-1}\|_{\mathrm{op}}
    \leq 1.
\]
Therefore,
\[
\begin{aligned}
    \left\|
        B(\beta)(\overline{T}(0)-iI)^{-1}
    \right\|_{\mathrm{op}}
    &\leq
    \|B(\beta)\|_{\mathrm{op}}
    \|(\overline{T}(0)-iI)^{-1}\|_{\mathrm{op}}\\
    &\leq
    C|\beta|^{q+2p}.
\end{aligned}
\]
Since $q+2p>0$, the right-hand side tends to zero as
$\beta\to0$. Hence there exists $\rho_0\in(0,\rho]$ such that
\[
    \left\|
        B(\beta)(\overline{T}(0)-iI)^{-1}
    \right\|_{\mathrm{op}}
    <1
\]
whenever $|\beta|<\rho_0$.

Recall the Neumann-series representation of the resolvent: if $T$ is
a bounded operator and $|\lambda|>\|T\|_{\mathrm{op}}$, then
\[
    (\lambda I-T)^{-1}
    =
    \frac{1}{\lambda}
    \sum_{n=0}^{\infty}
    \left(\frac{T}{\lambda}\right)^n,
\]
where the series converges in operator norm; see
\cite[equation~(VI.2), p. 191, Section~VI.3, Vol. 1 ]{RS}.
Taking $T=-K$ and $\lambda=1$, if
$\|K\|_{\mathrm{op}}<1$ then
\[
    (I+K)^{-1}
    =
    \sum_{n=0}^{\infty}(-K)^n.
\]
Applying this result with
\[
    K
    =
    B(\beta)(\overline{T}(0)-iI)^{-1}
\]
shows that
\[
    I+B(\beta)(\overline{T}(0)-iI)^{-1}
\]
is invertible for $|\beta|<\rho_0$.

Now,
\[
\begin{aligned}
    \overline{T}(\beta)-iI
    &=
    \overline{T}(0)-iI+B(\beta)\\
    &=
    \left[
        I+B(\beta)(\overline{T}(0)-iI)^{-1}
    \right]
    (\overline{T}(0)-iI),
\end{aligned}
\]
where the second equality follows from
\[
    B(\beta)(\overline{T}(0)-iI)^{-1}
    (\overline{T}(0)-iI)
    =
    B(\beta).
\]
Both factors in this factorisation are invertible, so 
$\overline{T}(\beta)-iI$ is invertible, with
\[
    (\overline{T}(\beta)-iI)^{-1}
    =
    (\overline{T}(0)-iI)^{-1}
    \left[
        I+B(\beta)(\overline{T}(0)-iI)^{-1}
    \right]^{-1}.
\]
As the factors on the right-hand side are bounded, the left-hand side inverse is also bounded. By definition we have
\[
    i\in\rho\bigl(\overline{T}(\beta)\bigr)
\]
for every $|\beta|<\rho_0$.
\end{proof}

\begin{lemma}\label{lemma_rescaled_analytic_family_A}
Let $s=p/q\in(0,1)$ be rational. Then the family
$\overline{T}(\beta)$ given by
\eqref{rescaled_family_pert_op_Taylor} is an analytic family of type~(A)
for $\beta$ in a sufficiently small complex neighbourhood of zero. For
real $\beta$ in this neighbourhood, it is self-adjoint.
\end{lemma}

\begin{proof}
Let $D:=\prod_{i=1}^k D\bigl((-\Delta)^s_{\mathrm{res},I}\bigr)$. After restricting
to a sufficiently small neighbourhood of zero, Lemma~\ref{rescaled_pert_op_nonempty_resolvent} shows that
$\overline{T}(\beta)$ is closed on the common domain $D$ and has non-empty resolvent set for every $\beta$ in this neighbourhood. In order to verify
Definition~\ref{def_analytic_family_A}, all that remains to be proved is
that for each $\fvect{u}\in D$, the $\mathcal H$-valued map
\[
    \beta\longmapsto\overline{T}(\beta)\fvect{u}
\]
is analytic for $\beta$ sufficiently close to zero.

Write
\[
    \overline{T}(\beta)
    =
    \overline{T}(0)+B(\beta),
\]
where
\begin{align*}
    [B(\beta)\fvect{u}]_i(z)
    &=
    \beta^{q+2p}V_i u_i(z) \\
    &\quad
    -c_s\beta^{q+2p}
    \sum_{j\neq i}
    \int_I
    K_{ij}(\beta;z,w)u_j(w)\,dw .
\end{align*}
For $i\neq j$, define
\[
    [C_{ij}(\beta)u](z)
    :=
    \int_I K_{ij}(\beta;z,w)u(w)\,dw .
\]
By \eqref{definition_analytic_kernel},
\[
    K_{ij}(\beta;z,w)
    =
    \frac{1}{|d_{ij}|^{1+2s}}
    \sum_{n=0}^{\infty}
    \binom{-(1+2s)}{n}
    \left(
        \frac{z-w}{d_{ij}}
    \right)^n
    \beta^{qn}.
\]
For each $n\geq 0$, define the integral operator
$C_{ij,n}:L^2(I)\to L^2(I)$ by
\[
    [C_{ij,n}u](z)
    :=
    \frac{1}{|d_{ij}|^{1+2s}}
    \binom{-(1+2s)}{n}
    \int_I
    \left(
        \frac{z-w}{d_{ij}}
    \right)^n
    u(w)\,dw.
\]
Thus $C_{ij,n}$ is the integral operator corresponding to the
coefficient of $\beta^{qn}$ in the power-series expansion of
$K_{ij}(\beta;z,w)$. This leads to the operator-valued power series
\[
    C_{ij}(\beta)
    =
    \sum_{n=0}^{\infty}
    \beta^{qn}C_{ij,n},
\]
whose convergence in operator norm we now verify.

Recall from the proof of
Lemma~\ref{rescaled_pert_op_nonempty_resolvent} that
\[
    d_*:=\min_{i\neq j}|d_{ij}|>0
\]
and that $\rho>0$ was chosen sufficiently small that
\[
    \theta:=\frac{2\rho^q}{d_*}<1.
\]
Since $I=(-1,1)$, we have $|z-w|\leq 2$. Hence, for
$|\beta|<\rho$,
\[
    \left|
        \frac{\beta^q(z-w)}{d_{ij}}
    \right|
    \leq
    \theta<1.
\]
Consequently, the kernel of the $n$th term
$\beta^{qn}C_{ij,n}$ satisfies
\[
\begin{aligned}
    &\frac{|\beta|^{qn}}{|d_{ij}|^{1+2s}}
    \left|
        \binom{-(1+2s)}{n}
    \right|
    \left|
        \frac{z-w}{d_{ij}}
    \right|^n\\
    &\qquad\leq
    \frac{1}{d_*^{1+2s}}
    \left|
        \binom{-(1+2s)}{n}
    \right|
    \theta^n
    =:M_n .
\end{aligned}
\]
As shown in the proof of
Lemma~\ref{rescaled_pert_op_nonempty_resolvent},
\[
    \sum_{n=0}^{\infty}M_n<\infty.
\]

The same integral-operator estimate used there gives
\[
    \|\beta^{qn}C_{ij,n}u\|_{L^2(I)}
    \leq
    |I|M_n\|u\|_{L^2(I)}
\]
for every $u\in L^2(I)$. Taking the supremum over
$\|u\|_{L^2(I)}=1$ therefore yields
\[
    \|\beta^{qn}C_{ij,n}\|_{\mathrm{op}}
    \leq
    |I|M_n.
\]
It follows that
\[
    \sum_{n=0}^{\infty}
    \|\beta^{qn}C_{ij,n}\|_{\mathrm{op}}
    \leq
    |I|
    \sum_{n=0}^{\infty}M_n
    <\infty.
\]
Hence, 
\[\sum_{n=0}^{\infty}\beta^{qn}C_{ij,n}\] 
converges absolutely in the operator norm of $\mathcal B(L^2(I))$. 
Its sum is the integral operator
$C_{ij}(\beta)$ associated with the kernel
$K_{ij}(\beta;z,w)$.
Therefore
\[
    C_{ij}(\beta)
    =
    \sum_{n=0}^{\infty}
    \beta^{qn}C_{ij,n}
\]
is an operator-norm convergent power series in $\beta$ for
$|\beta|<\rho$. Since a power series with coefficients in a Banach
space may be differentiated term by term within its radius of
convergence, it follows that
\[
    \beta\longmapsto C_{ij}(\beta)
\]
is analytic in operator norm, and hence analytic as a
$\mathcal B(L^2(I))$-valued function in the sense of Kato; see
\cite[Chapter~VII, \S1.1]{Kato}.

The potential term depends polynomially on $\beta$, and there are only
finitely many components and interaction terms. As each
$\beta\mapsto C_{ij}(\beta)$ is analytic in operator norm, it follows
that
\[
    \beta\longmapsto B(\beta)
\]
is analytic as a $\mathcal B(\mathcal H)$-valued function for
$|\beta|<\rho$. Thus $B(\beta)$ is a bounded-operator-valued analytic function. 

In particular, the derivative $B'(\beta)$ exists in operator norm, so
\[
    \left\|
        \frac{B(\beta+h)-B(\beta)}{h}
        -
        B'(\beta)
    \right\|_{\mathrm{op}}
    \longrightarrow 0
    \qquad\text{as }h\to0.
\]
Consequently, for every fixed $\fvect{u}\in\mathcal H$,
\[
\begin{aligned}
    &
    \left\|
        \frac{B(\beta+h)\fvect{u}-B(\beta)\fvect{u}}{h}
        -
        B'(\beta)\fvect{u}
    \right\|_{\mathcal H}
    \\
    &\qquad\leq
    \left\|
        \frac{B(\beta+h)-B(\beta)}{h}
        -
        B'(\beta)
    \right\|_{\mathrm{op}}
    \|\fvect{u}\|_{\mathcal H}
    \longrightarrow 0.
\end{aligned}
\]
Thus, for every fixed $\fvect{u}\in\mathcal H$,
\[
    \beta\longmapsto B(\beta)\fvect{u}
\]
is an $\mathcal H$-valued analytic function for $|\beta|<\rho$.

We now return to the generally unbounded family
\[
    \overline{T}(\beta)
    =
    \overline{T}(0)+B(\beta),
\]
whose common operator domain is $D$. For each
$\fvect{u}\in D$, the vector $\overline{T}(0)\fvect{u}$ is independent
of $\beta$, while $\beta\mapsto B(\beta)\fvect{u}$ is analytic.
Therefore
\[
    \beta\longmapsto
    \overline{T}(\beta)\fvect{u}
    =
    \overline{T}(0)\fvect{u}
    +
    B(\beta)\fvect{u}
\]
is an $\mathcal H$-valued analytic function of $\beta$ near zero for every fixed
$\fvect{u}\in D$.

Finally, let $\rho_0\leq\rho$ be the radius obtained in
Lemma~\ref{rescaled_pert_op_nonempty_resolvent}, so that for
$|\beta|<\rho_0$ the operators $\overline{T}(\beta)$ are closed on the
common domain $D$ and have non-empty resolvent set. On this smaller
complex neighbourhood, the preceding vector-valued analyticity
therefore verifies Definition~\ref{def_analytic_family_A}. Hence
\[
    \bigl\{\overline{T}(\beta):|\beta|<\rho_0\bigr\}
\]
is an analytic family of type~(A).

It remains to establish self-adjointness for real $\beta$. The
unperturbed operator $\overline{T}(0)=\bigoplus_{i=1}^k(-\Delta)^s_{\mathrm{res},I}$ is self-adjoint on $D$. Moreover, the factor $\beta^{q+2p}$ is real, the potentials $V_i$ are real-valued, and
\[
    K_{ij}(\beta;z,w)
    =
    K_{ji}(\beta;w,z).
\]
Hence the bounded operator $B(\beta)$ is symmetric on $\mathcal H$,
and therefore
\[
    \overline T(\beta)
    =
    \overline T(0)+B(\beta)
\]
is symmetric on $D$.

We verify self-adjointness directly. Let
\[
    \fvect{v}\in D\bigl(\overline T(\beta)^*\bigr).
\]
Then, for every $\fvect{u}\in D$,
\[
\begin{aligned}
    \langle \overline T(0)\fvect{u},\fvect{v}\rangle_{\mathcal H}
    &=
    \langle \overline T(\beta)\fvect{u},\fvect{v}\rangle_{\mathcal H}
    -
    \langle B(\beta)\fvect{u},\fvect{v}\rangle_{\mathcal H}\\
    &=
    \langle \fvect{u},
        \overline T(\beta)^*\fvect{v}\rangle_{\mathcal H}
    -
    \langle \fvect{u},B(\beta)\fvect{v}\rangle_{\mathcal H}\\
    &=
    \langle \fvect{u},
        \overline T(\beta)^*\fvect{v}
        -
        B(\beta)\fvect{v}\rangle_{\mathcal H}.
\end{aligned}
\]
By the definition of the adjoint,
\[
    \fvect{v}\in D\bigl(\overline T(0)^*\bigr).
\]
Since $\overline T(0)$ is self-adjoint,
\[
    D\bigl(\overline T(0)^*\bigr)
    =
    D\bigl(\overline T(0)\bigr)
    =
    D.
\]
It follows that
\[
    D\bigl(\overline T(\beta)^*\bigr)\subset D.
\]
The reverse inclusion follows from the symmetry of
$\overline T(\beta)$, and therefore
\[
    D\bigl(\overline T(\beta)^*\bigr)=D.
\]
Thus $\overline T(\beta)$ is self-adjoint for every real $\beta$ in the neighbourhood under consideration.
\end{proof}

\begin{lemma}\label{limiting_rescaled_ground_state_degenerate}
The limiting re-scaled operator is
\begin{align}\label{limiting_rescaled_operator}
    \overline{T}(0)
    =\bigoplus_{i=1}^k(-\Delta)^s_{\mathrm{res},I}.
\end{align}
It acts on $\prod_{i=1}^kL^2(I)$ with operator domain
\[
    D(\overline{T}(0))
    =\prod_{i=1}^kD\!\left((-\Delta)^s_{\mathrm{res},I}\right).
\]
If $\lambda_1$ is the ground state eigenvalue of
$(-\Delta)^s_{\mathrm{res},I}$, then $\lambda_1$ is a discrete
eigenvalue of $\overline{T}(0)$ of multiplicity exactly $k$.
\end{lemma}

\begin{proof}
Setting $\beta=0$ in
\eqref{rescaled_family_pert_op_Taylor} eliminates both the potential
and the interactions between distinct components, giving
\eqref{limiting_rescaled_operator}.

Let $\phi_1$ be the positive $L^2(I)$-normalised ground state
eigenfunction. Since the ground state eigenvalue of each one-interval
operator is simple,
\[
    \ker(\overline{T}(0)-\lambda_1)
    =
    \operatorname{span}
    \left\{
        \fvect{e}_1\phi_1,\ldots,\fvect{e}_k\phi_1
    \right\},
\]
where $\fvect{e}_i\phi_1$ has $\phi_1$ in the $i$th component and zero
in every other component. This eigenspace has dimension $k$.

Moreover, the next eigenvalue of the direct-sum operator is the second
one-interval eigenvalue $\lambda_2>\lambda_1$ (see, for example, \cite{BK}), so $\lambda_1$ is
discrete and isolated.
\end{proof}

We are now in a position to use the Kato--Rellich theory about the degenerate ground state of the limiting re-scaled operator.
\begin{theorem}\label{rescaled_op_existence_pert_evalues}
Let $s=p/q\in(0,1)$ be rational. Then, for $\beta$ sufficiently close to zero, the re-scaled family $\overline{T}(\beta)$ admits $k$ eigenvalue branches, counted with multiplicity, that are analytic in $\beta$ and satisfy
\[
    \lambda^{(j)}(0)=\lambda_1,
    \qquad j=1,\ldots,k.
\]
These are the only eigenvalues of $\overline{T}(\beta)$ in a sufficiently small neighbourhood of $\lambda_1$.
\end{theorem}
\begin{proof}
By Lemma \ref{lemma_rescaled_analytic_family_A}, $\overline{T}(\beta)$ is an analytic family of type (A), and therefore an analytic family in the sense of Kato. The same lemma proves that it is self-adjoint for real $\beta$ sufficiently close to zero. Finally, Lemma \ref{limiting_rescaled_ground_state_degenerate} shows that $\lambda_1$ is a discrete eigenvalue of multiplicity $k$ of $\overline{T}(0)$. The result follows from Theorem \ref{Kato_Rellich_degenerate_evalue}.
\end{proof}

\subsection{Constructing the infinite potential well counterexample}\label{reg_pert_calc_s_half}
Theorem \ref{rescaled_op_existence_pert_evalues} enables us to apply
Rayleigh--Schr\"{o}dinger regular perturbation theory to the $k$
eigenvalue branches emerging from the degenerate ground state of
$\overline{T}(0)$. We now take $s=1/2$ and $k=3$. Under the parametrisation
$s=p/q$ and $\beta=\epsilon^{1/q}$ introduced above, this corresponds to $p=1$ and $q=2$, giving $\beta=\sqrt{\epsilon}$. 
For notational simplicity, write $\lambda(\beta)$ for an eigenvalue
of the re-scaled operator $\overline{T}(\beta)$. It corresponds to
$\epsilon\lambda_\epsilon=\beta^2\lambda_\epsilon$ in the original
problem. The exact re-scaled problem is
\begin{align}\label{rescaled_evalue_problem_k_3}
 \lambda(\beta)u_i^\beta(z)
 &=(-\Delta)^{1/2}_{\mathrm{res},I}u_i^\beta(z)
   +\beta^4V_iu_i^\beta(z)\nonumber\\
 &\quad-c\beta^4\sum_{j\neq i}\int_I
 K_{ij}(\beta;z,w)u_j^\beta(w)\,dw,
\end{align}
where $i=1,2,3$ and
\[
 K_{ij}(\beta;z,w)
 =\frac{1}{|x_i-x_j|^2}
  \left(1+\frac{\beta^2(z-w)}{x_i-x_j}\right)^{-2}.
\]
Expanding the interaction kernel near $\beta=0$ gives
\[
    K_{ij}(\beta;z,w)
    =\frac{1}{|x_i-x_j|^2}+O(\beta^2),
\]
uniformly for $z,w\in I$. Hence the corresponding interaction
operator satisfies
\[
    C_{ij}(\beta)=C_{ij}(0)+O(\beta^2)
\]
in operator norm, and therefore
\[
    \beta^4 C_{ij}(\beta)
    =
    \beta^4 C_{ij}(0)+O(\beta^6).
\]
Since the potential term also begins at order $\beta^4$, it follows
that the re-scaled operator has the expansion
\[
    \overline{T}(\beta)
    =
    \overline{T}(0)+\beta^4T_4+O(\beta^6)
\]
for a bounded operator $T_4$. Thus the first non-trivial perturbation
of $\overline{T}(0)$ occurs at order $\beta^4$.

To carry out the perturbation calculation, we use the above expansion
of $\overline{T}(\beta)$ together with power-series expansions of its
eigenvalue $\lambda(\beta)$ and the corresponding eigenfunction
$\fvect{u}^{\beta}$:
\begin{align*}
    \lambda(\beta)
      &=\lambda^0+\beta\lambda^1+\beta^2\lambda^2
        +\beta^3\lambda^3+\beta^4\lambda^4+\cdots,\\
    \fvect{u}^\beta
      &=\fvect{u}^0+\beta\fvect{u}^1
        +\beta^2\fvect{u}^2+\beta^3\fvect{u}^3
        +\beta^4\fvect{u}^4+\cdots.
\end{align*}
At zeroth order, the unperturbed problem is
\[
    \overline{T}(0)\fvect{u}^0
    =
    \lambda^0\fvect{u}^0.
\]
As shown in Lemma~\ref{limiting_rescaled_ground_state_degenerate},
$\lambda^0$ is the threefold degenerate ground-state eigenvalue of
$\overline{T}(0)$, with eigenspace
\[
    \ker\bigl(\overline{T}(0)-\lambda^0I\bigr)
    =
    \operatorname{span}\left\{
        \fvect{e}_1\phi_1,\,
        \fvect{e}_2\phi_1,\,
        \fvect{e}_3\phi_1
    \right\},
\]
and $\phi_1$ is the $L^2(I)$-normalised ground-state eigenfunction of
$(-\Delta)^{1/2}_{\mathrm{res},I}$. Hence, the most general zeroth-order
eigenfunction has the form
\[
    \fvect{u}^0
    =
    \sum_{i=1}^3\alpha^i\fvect{e}_i\phi_1
    =
    \begin{bmatrix}
        \alpha^1\phi_1\\
        \alpha^2\phi_1\\
        \alpha^3\phi_1
    \end{bmatrix}.
\]
We now substitute the expansions of $\overline{T}(\beta)$,
$\lambda(\beta)$, and $\fvect{u}^{\beta}$ into
\[
    \overline{T}(\beta)\fvect{u}^{\beta}
    =
    \lambda(\beta)\fvect{u}^{\beta}
\]
and compare coefficients of like powers of $\beta$. This gives an
equation at each order in the perturbation expansion, which is then used
to determine the leading non-trivial corrections.

Since the operator perturbation in \eqref{rescaled_evalue_problem_k_3}
begins at order $\beta^4$, comparing the coefficients at orders
$\beta$, $\beta^2$, and $\beta^3$, together with the corresponding
solvability conditions, gives
\[
    \lambda^1=\lambda^2=\lambda^3=0.
\]
At order $\beta^4$, the $i$th component satisfies
\begin{align}\label{fourth_order_pert}
 \left((-\Delta)^{1/2}_{\mathrm{res},I}-\lambda^0\right)u_i^4(z)
 &=\lambda^4\alpha^i\phi_1(z)-V_i\alpha^i\phi_1(z)\nonumber\\
 &\quad+c\sum_{\ell\neq i}\frac{\alpha^\ell}{|x_i-x_\ell|^2}
       \int_I\phi_1(w)\,dw.
\end{align}
Taking the $L^2(I)$ inner product with $\phi_1$ and using
self-adjointness yields the solvability condition
\begin{align}\label{effective_matrix_equations}
    \lambda^4\alpha^i
    =V_i\alpha^i
     -c\mu_1^2\sum_{\ell\neq i}
       \frac{\alpha^\ell}{|x_i-x_\ell|^2},
    \qquad
    \mu_1:=\int_I\phi_1(z)\,dz>0.
\end{align}
Thus the solvability condition reduces to the three-by-three matrix
eigenvalue problem
\begin{align}\label{evalue_eqn_pert_matrix}
    \widehat{\M}\fvect{\alpha}
    =\lambda^4\fvect{\alpha},
\end{align}
where
\[
    \fvect{\alpha}
    =
    (\alpha^1,\alpha^2,\alpha^3)^T,
\]
and
\begin{align}\label{correct_effective_matrix}
 \widehat{\M}
 =\begin{bmatrix}
 V_1 & -\dfrac{c\mu_1^2}{|x_1-x_2|^2}
     & -\dfrac{c\mu_1^2}{|x_1-x_3|^2}\\[8pt]
 -\dfrac{c\mu_1^2}{|x_2-x_1|^2} & V_2
     & -\dfrac{c\mu_1^2}{|x_2-x_3|^2}\\[8pt]
 -\dfrac{c\mu_1^2}{|x_3-x_1|^2}
     & -\dfrac{c\mu_1^2}{|x_3-x_2|^2} & V_3
 \end{bmatrix}.
\end{align}

The matrix \eqref{correct_effective_matrix} has the form analysed in Section \ref{putative_counterexample}. Its off-diagonal entries are negative and, since $x_1<x_2<x_3$,
\[
    |\widehat M_{13}|<|\widehat M_{12}|,
    \qquad
    |\widehat M_{13}|<|\widehat M_{23}|.
\]
Choosing the potential values so that $V_2<V_1,V_3$ therefore gives three simple eigenvalues of $\widehat{\M}$. The coefficient vector associated with the smallest eigenvalue has constant sign, the coefficient vector associated with the second eigenvalue has signs $(-,+,-)$ or $(+,-,+)$, and the coefficient vector associated with the largest eigenvalue has one sign change.

It follows that, for $0<\beta\ll1$, the three eigenvalue branches
emerging from the $k=3$ degenerate ground-state eigenvalue satisfy
\[
    \lambda_j(\beta)
    =
    \lambda^0+\beta^4\lambda_j^4+o(\beta^4),
    \qquad j=1,2,3,
\]
where $\lambda_1^4,\lambda_2^4,\lambda_3^4$ are the distinct
eigenvalues of the effective matrix $\widehat{\M}$. Hence the
threefold degeneracy is split, and the corresponding eigenvalue
branches are simple for all sufficiently small $\beta>0$.

If
\[
    \boldsymbol{\alpha}_j
    =
    \bigl(
        \alpha_j^1,
        \alpha_j^2,
        \alpha_j^3
    \bigr)^T
\]
is an eigenvector of $\widehat{\M}$ corresponding to
$\lambda_j^4$, then it determines the zeroth-order term of the
associated re-scaled eigenfunction branch; namely,
\[
    \fvect{u}_j^0
    =
    \sum_{\ell=1}^3
    \alpha_j^\ell\fvect{e}_\ell\phi_1,
    \qquad j=1,2,3.
\]
Consequently,
\[
    \fvect{u}_j(\beta)
    \longrightarrow
    \fvect{u}_j^0
    =
    \sum_{\ell=1}^3
    \alpha_j^\ell\fvect{e}_\ell\phi_1
    =
    \begin{bmatrix}
        \alpha_j^1\phi_1\\
        \alpha_j^2\phi_1\\
        \alpha_j^3\phi_1
    \end{bmatrix}
    \qquad
    \text{as }\beta\to0^+,
\]
since all positive-order terms in the expansion of
$\fvect{u}_j(\beta)$ vanish as $\beta\to0^+$. Thus the eigenvectors of
$\widehat{\M}$ determine the leading-order linear combinations of the
degenerate orthonormal basis $\mathcal B$. Since $\phi_1$ is strictly positive on $I$, the sign of each component is determined by the sign of the corresponding entry of $\boldsymbol{\alpha}_j$. In particular, the limiting combination associated with the lowest branch can be chosen strictly positive across the three components, that associated with the highest branch has one sign change, while that associated with the second branch has two sign changes. The latter is the sign pattern required for the infinite potential well counterexample.

The same construction of the analytic re-scaled family applies to every rational $s\in(0,1)$. For general rational $s=p/q$, the first non-trivial perturbation occurs at order $\beta^{q+2p}$ and produces the analogous effective matrix with off-diagonal entries proportional to $-|x_i-x_j|^{-(1+2s)}$. We restrict to $s=1/2$ and $k=3$, since this gives the simplest instance of the construction in which the perturbation produces the sign pattern required for the counterexample.

\section{Finite potential well counterexample}\label{finite_well_counterexample}

The conclusion of the previous section shows that, for $\epsilon>0$ sufficiently small, the first three eigenvalues of the re-scaled infinite potential well problem are simple and the corresponding eigenfunctions are governed, to first non-trivial order, by the finite-dimensional matrix analysed in Section \ref{putative_counterexample}. Since the re-scaling is unitary, it preserves the signs of the eigenfunctions on the corresponding intervals. In particular, the second eigenfunction of the original infinite potential well problem has two sign changes when the potential values are non-convex and satisfy $V_2<V_1,V_3$.

\begin{figure}
    \centering
    \includegraphics[width=12cm]{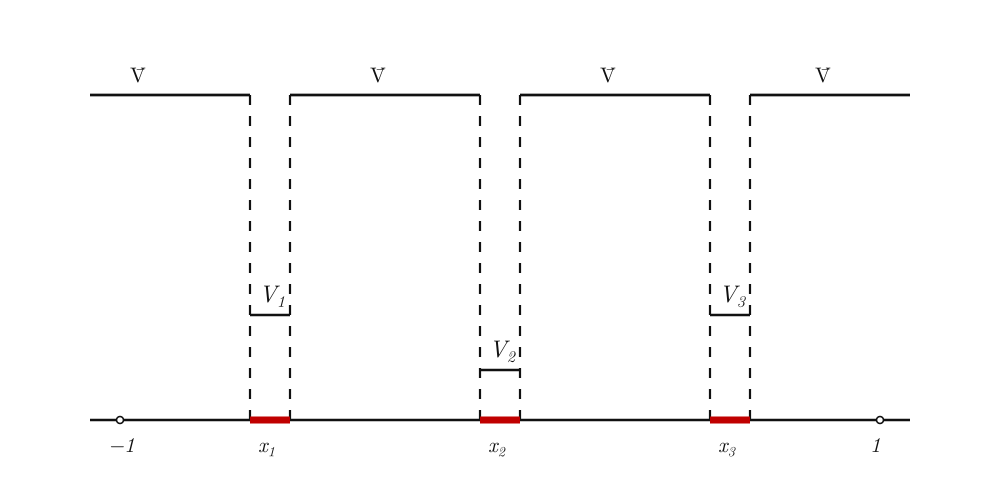}
    \caption{A triple finite potential well with constant potential value $\bar{V}=1/\delta\gg V_1,V_2,V_3$ off the subintervals, which are coloured in red.} 
    \label{fig_finite_potential}
\end{figure}

To complete the counterexample, we now fix $\epsilon$ sufficiently small so that the subintervals $I_i$ have uniform fixed widths, and `perturb' the infinite potential well $V^\infty$ such that it takes a very large and positive value off the finite union of mutually disjoint subintervals in $(-1,1)$. See Figure \ref{fig_finite_potential}. Our aim is to show that as the finite potential well depth grows to infinity, the eigenfunctions of the resulting `perturbed' finite potential well eigenvalue problem converge in $L^2$ to the eigenfunctions of the `unperturbed' infinite potential well eigenvalue problem. The perturbed eigenfunctions can thus be made arbitrarily close in $L^2$ to the unperturbed eigenfunctions in this limit. In particular, the perturbed second eigenfunction inherits the sign change properties of the unperturbed second eigenfunction and so has more than one sign change.\\

Precisely, let $U$ denote a finite union of mutually disjoint open subintervals $I_i$, for $i=1,\ldots,N$, in the bounded interval $I:=(-1,1)$. The subintervals have uniform fixed widths that are very small compared to the distances between adjacent subintervals. For any $s$ in $(0,1)$, let $u_i$ denote an eigenfunction of the $i$th eigenvalue $\lambda_i$ of the infinite potential well eigenvalue problem

\begin{equation}\label{evalue_problem_infinite_well}
 \left\{\begin{aligned}
        ((-\Delta)^s+V(x))u_i(x)&=\lambda_i u_i(x), &x\in U\\
        u_i(x)&=0, &x\in\R\setminus U,
       \end{aligned}
 \right.
\end{equation}
where $V:U\to\R$ denotes a non-negative bounded (and smooth) potential on $U$. Recall from Section \ref{ch_4_background} that the eigenvalues $\lambda_i$ can be degenerate and form an increasing sequence $\lambda_i\leq\lambda_{i+1}$ in the index $i$. 
Moreover, the corresponding eigenfunctions can be chosen to be orthonormal in $L^2(U)$. 

For each $0<\delta\ll 1$, let $u_{i,\delta}$ be an eigenfunction of the $i$th eigenvalue $\lambda_{i,\delta}$ of the finite potential well eigenvalue problem
\begin{equation}\label{evalue_problem_finite_well}
 \left\{\begin{aligned}
        ((-\Delta)^s+V(x))u_{i,\delta}(x)&=\lambda_{i,\delta} u_{i,\delta}(x), &x\in U\\
        \left((-\Delta)^s+\frac{1}{\delta}\right)u_{i,\delta}(x)&=\lambda_{i,\delta} u_{i,\delta}(x), &x\in I\setminus U\\
        u_{i,\delta}(x)&=0, &x\in\R\setminus I.
       \end{aligned}
 \right.
\end{equation}
The eigenvalues $\lambda_{i,\delta}$ are also 
increasing in $i$. The corresponding eigenfunctions $u_{i,\delta}$ can be chosen to be orthonormal in $L^2(I)$.

Finally, recalling the Hilbert space structure of $H_0^s(I)$ induced by the Dirichlet energy $\mathcal E_{\mathrm{res}}$ from Section~\ref{subs_func_an}, we first note that
$H_0^s(U)$ is a closed subspace of $H_0^s(I)$. Indeed,
$H_0^s(U)\subset H_0^s(I)$, and if
$v_n\in H_0^s(U)$ with
\[
    v_n\longrightarrow v
    \qquad\text{in }H_0^s(I),
\]
then Proposition~\ref{prop_equiv_norms} implies that
$v_n\to v$ in $H^s(\R)$, and therefore also in $L^2(\R)$. Since each
$v_n$ vanishes almost everywhere on $\R\setminus U$,
\[
    \|v\|_{L^2(\R\setminus U)}
    =
    \|v-v_n\|_{L^2(\R\setminus U)}
    \leq
    \|v-v_n\|_{L^2(\R)}
    \longrightarrow0.
\]
Thus $v=0$ almost everywhere on $\R\setminus U$, so
$v\in H_0^s(U)$. Therefore $H_0^s(U)$ is closed in $H_0^s(I)$.

It follows that
\[
    H_0^s(I)
    =
    H_0^s(U)\oplus H_0^s(U)^{\perp_{\mathcal E}},
\]
where the orthogonal complement is taken with respect to the inner
product $\mathcal E_{\mathrm{res}}$; see, for example, \cite[Sec.~3.3, Thm.~3.3-4]{Krey}.

For each $u_{i,\delta}\in H_0^s(I)$, we may accordingly write
\[
    u_{i,\delta}
    =
    v_{i,\delta}+w_{i,\delta},
\]
where
\[
    v_{i,\delta}\in H_0^s(U),
    \qquad
    w_{i,\delta}\in H_0^s(U)^{\perp_{\mathcal E}}.
\]
Thus
\[
    \mathcal E_{\mathrm{res}}[w_{i,\delta},\varphi]=0
    \qquad
    \text{for every }\varphi\in H_0^s(U).
\]
In particular,
\[
    \mathcal E_{\mathrm{res}}[v_{i,\delta},w_{i,\delta}]=0.
\]
Since $v_{i,\delta}=0$ on $\R\setminus U$, we also have
$w_{i,\delta}=u_{i,\delta}$ there. We see that the preceding orthogonality condition is the weak formulation of
\begin{equation}\label{BVP_w}
 \left\{
 \begin{aligned}
        (-\Delta)^s w_{i,\delta}(x)&=0, &x\in U,\\
        w_{i,\delta}(x)&=u_{i,\delta}(x), &x\in\R\setminus U.
 \end{aligned}
 \right.
\end{equation}
Existence of solutions to BVP (\ref{BVP_w}) is well-developed; see, for instance, \cite[Sec. 3]{RO}. 

Central to showing $L^2$ convergence of the finite well eigenfunctions $u_{i,\delta}$ to the infinite well eigenfunctions $u_i$ is an estimate that bounds the $L^2$ norm of $w_{i,\delta}$ on $U$ by a higher $L^p$ norm on $\R\setminus\overline{U}$ for $p$ depending on $s$.
\begin{lemma}\label{L2_estimate_w} Let $w_{i,\delta}$ be as defined above and satisfy BVP (\ref{BVP_w}). Then the following estimate 
\begin{align}\label{L2_bound_w}
\|w_{i,\delta}\|_{L^2(U)}\leq \tilde{C}\|w_{i,\delta}\|_{L^p(\R\setminus\overline{U})}
\end{align}
holds for $p> \frac{1}{1-s}$ and $s>1/2$, and $p>2$ for $s\leq 1/2$. $\tilde{C}$ is a large positive constant depending on $s$ and $U$.
\end{lemma}
\begin{proof}
    To simplify notation, set $w:=w_{i,\delta}$ and $u:=u_{i,\delta}$. For any $y$ in $\R$, let $\delta(y):=dist(y, \partial U)=\inf\{|y-x|: x\in\partial U\}$. 
    
    From \cite[Th. 1.2.3]{Ab1}, we know that the (unique) pointwise solution $w$ of BVP (\ref{BVP_w}) has a representation in terms of the $s$-harmonic extension of the boundary condition $u$ into $U$ given by
    \begin{align}\label{rep_w}
        w(x)=\int_{\R\setminus\overline{U}}\Gamma_s(x,y)u(y)dy=\int_{I\setminus\overline{U}}\Gamma_s(x,y)u(y)dy,
    \end{align}
    where the second equality follows from $u=0$ on $\R\setminus I$. Moreover, from \cite[Rm. 1.2.5]{Ab1}, the integral kernel $\Gamma_s(\cdot,\cdot)$ is bounded, and the upper bound can be simplified for the ensuing calculations: for $x\in U$ and $y\in \R\setminus\overline{U}$,
    \begin{align}\label{ineq_Gamma_s}
        \Gamma_s(x,y)&\leq \frac{C\delta(x)^s}{\delta(y)^s(1+\delta(y))^s|x-y|}\nonumber\\
        &\leq \frac{C\delta(x)^s}{\delta(y)^s|x-y|},
    \end{align}
    and the positive constant $C$ depends on $s$ and $U$.
    From any point $x$ in $U$ we must pass through at least one boundary point to reach a point $y$ in $\R\setminus\overline{U}$, so we also have the following lower bounds on the distance between $x$ and $y$:
    \begin{align}\label{lower_bound_xy}
        |x-y|\geq |\delta(y)+\delta(x)|\geq \delta(x),\delta(y).
    \end{align}
    Using equation (\ref{rep_w}), inequality (\ref{ineq_Gamma_s}), Fubini's theorem, and H\"{o}lder's inequality, we can calculate an upper bound for $\|w\|_{L^2(U)}^2$ given by
    \begin{align}\label{L2_inquality_w}
        \int_{U}|w(x)|^2dx&\leq\int_{I\setminus\overline{U}}\int_{I\setminus\overline{U}}Q_s(y,z)u(y)u(z)dydz,\qquad x\in U\nonumber\\
        &\leq\left(\int_{I\setminus\overline{U}}\int_{I\setminus\overline{U}}Q_s(y,z)^qdydz\right)^{1/q}\left(\int_{I\setminus\overline{U}}\int_{I\setminus\overline{U}}|u(y)|^p|u(z)|^pdydz\right)^{1/p}\nonumber\\
        &\leq\left(\int_{I\setminus\overline{U}}\int_{I\setminus\overline{U}}Q_s(y,z)^qdydz\right)^{1/q}\left(\int_{I\setminus\overline{U}}|u(y)|^pdy\right)^{2/p}\nonumber\nonumber\\
        &=\left(\int_{I\setminus\overline{U}}\int_{I\setminus\overline{U}}Q_s(y,z)^qdydz\right)^{1/q}\|u\|^2_{L^p(\R\setminus\overline{U})}
    \end{align}
    with integral kernel 
    \begin{align}\label{kernal_Qs}
        Q_s(y,z)&=\frac{C^2}{\delta(y)^s\delta(z)^s}\int_{U}\frac{\delta(x)^{2s}}{|x-y||x-z|}dx, \qquad y,z\in I\setminus\overline{U}.
    \end{align}    
    To complete the $L^2$ bound on $w$, we need to work out when $Q_s(\cdot,\cdot)^q$ is integrable for different values of $s$ and $q$.
    
    For $s>1/2$, we calculate
    \begin{align}\label{Q_bound_s>1/2}
        Q_s(y,z)&=\frac{C^2}{\delta(y)^s\delta(z)^s}\int_{U}\frac{\delta(x)^{2s}}{|x-y||x-z|}dx, \qquad y,z\in I\setminus\overline{U}\nonumber\\
        &\leq\frac{C^2}{\delta(y)^s\delta(z)^s}\int_{U}\frac{\delta(x)^{2s}}{\delta(x)\delta(x)}dx\nonumber\\
        &=\frac{C^2}{\delta(y)^s\delta(z)^s}\int_{U}\delta(x)^{2s-2}dx\nonumber\\
        &=\frac{C^2}{\delta(y)^s\delta(z)^s}\sum_{i=1}^N\int_{I_i}\delta(x)^{2s-2}dx\nonumber\\
        &\leq \frac{\tilde{C}^2}{\delta(y)^s\delta(z)^s},
    \end{align}
    where inequality (\ref{lower_bound_xy}) was applied in the second step. Since $\delta(x)$ approaches zero linearly as $x$ approaches an endpoint, $\delta(x)^{2s-2}$ is integrable on $U$, i.e., $s>1/2$, so $2s-2>-1$. Therefore,
    \begin{align}\label{s>1/2_Qq_bound}
        \int_{I\setminus\overline{U}}\int_{\times I\setminus\overline{U}} Q_s(y,z)^qdydz&\leq\int_{I\setminus\overline{U}}\int_{I\setminus\overline{U}}\frac{\tilde{C}^{2q}}{\delta(y)^{sq}\delta(z)^{sq}}dydz\nonumber\\
        &=\left(\int_{I\setminus\overline{U}}\frac{\tilde{C}^q}{\delta(y)^{sq}}dy\right)^2\nonumber\\
        &< \infty,
    \end{align}
    as the final integral is integrable for $q<1/s$.

    For $s\leq 1/2$, assume $R>0$ is such that every pair of boundary points is separated by a distance of at least $2R$, including those for the same subinterval. 
    There are four cases to consider:\\

    \noindent \textbf{Case 1:} $\delta(y),\delta(z)\geq R$. In this case, for any $x$ in $U$, $|z-x|\geq R$ and $|y-x|\geq R$, so we can calculate
    \begin{align}\label{case1_Qbound_s<1/2}
    Q_s(y,z)&=\frac{C^2}{\delta(y)^s\delta(z)^s}\int_{U}\frac{\delta(x)^{2s}}{|x-y||x-z|}dx, \qquad y,z\in I\setminus\overline{U}\nonumber\\
    &\leq\frac{C^2}{\delta(y)^s\delta(z)^s}\frac{1}{R^2}\int_{U}\delta(x)^{2s}dx\nonumber\\
    &\leq \frac{\tilde{C}^2}{\delta(y)^s\delta(z)^2},
    \end{align}
    as the final integral is integrable. Thus $Q_s(\cdot,\cdot)^q$ is integable for $q<1/s$ like for the case $s>1/2$.\\
    
    \noindent\textbf{Cases 2 and 3}: $\delta(y)<R$ and $\delta(z)\geq R$ (or vice versa). Here for any $x$ in $U$, $|z-x|\geq R$, while $|y-x|\geq \delta(x)$, therefore we get
    \begin{align}\label{case2_Qbound_s<1/2}
        Q_s(y,z)&=\frac{C^2}{\delta(y)^s\delta(z)^s}\int_{U}\frac{\delta(x)^{2s}}{|x-y||x-z|}dx, \qquad y,z\in I\setminus\overline{U}\nonumber\\
        &\leq\frac{C^2}{\delta(y)^s\delta(z)^s}\int_{U}\frac{\delta(x)^{2s}}{\delta(x)R}dx\nonumber\nonumber\\
        &\leq\frac{C^2}{\delta(y)^s\delta(z)^s}\frac{1}{R}\int_{U}\delta(x)^{2s-1}dx\nonumber\\
        &\leq \frac{\tilde{C}^2}{\delta(y)^s\delta(z)^2},
    \end{align}
    as $2s-1>-1$, so the final integral is integrable. Inequality (\ref{lower_bound_xy}) was also applied in the second line. It follows that $Q(\cdot,\cdot)^q$ is integrable for $q<1/s$ as before.\\

    \noindent\textbf{Case 4:} $\delta(y),\delta(z)<R$. If $y$ and $z$ are close to distinct boundary points, the situation reduces to one of the above cases. To see this, suppose $y$ and $z$ are respectively close to distinct boundary points $p_1$ and $p_2$ belonging to the same subinterval $I_i$ in $U$. Then for all $x$ in $I_i$ closer to $y$ than to $z$, we have $|z-x|=|z-p_2|+|p_2-x|\geq|z-p_2|+R\geq R$, while $|y-x|\geq\delta(x)$ (or vice versa). For all remaining points $x$ in $U\setminus I_i$, both $|z-x|\geq R$ and $|y-x|\geq R$. A similar argument applies if $y$ and $z$ are close to boundary points of distinct subintervals. So computing upper bounds for $Q_s(\cdot,\cdot)$ and $Q_s(\cdot,\cdot)^q$ proceeds in the same fashion as previously and yields the same result.
    
    The one situation remaining is when both $y$ and $z$ are close to the same boundary point $p$. In this case, we can partition $U$ into disjoint subsets $W=B_{R}(p)\cap U$ and its complement $V=U\setminus W$, where $B_R(p)$ denotes the open ball of radius $R$ centred at $p$. 
    Then for all $x$ belonging to the complement $V$ of $W$, $|y-x|=|y-p|+|p-x|\geq |y-p|+R> R$, and similarly, $|z-x|>R$, so we end up with the same $Q_s(\cdot,\cdot)$ and $Q_s(\cdot,\cdot)^q$ bounds as in Case 1. Therefore we only need to consider what happens on $W$.

    Without loss of generality (w.l.o.g), assume $p$ is the left-hand side end point of a subinterval $I_i$ in $U$. Then $W$ can be further partitioned into the subintervals $(p,p+M)$ and $[p+M,p+R)$, where $M=\max\{|y-p|,|z-p|\}$. We have the following inequalities:
    \begin{align*}
        |y-x||z-x|&\geq M|x-p| \qquad\qquad&\mbox{ on }(p,p+M)\\
        |y-x||z-x|&\geq |x-p|^2 \qquad\qquad&\mbox{ on }[p+M,p+R).
    \end{align*}
    The first inequality follows from the fact that one of $|y-p|$ or $|z-p|$ takes the maximum value in $M$, so one of $|y-x|$ or $|y-z|$ is bounded below by M, while both $|y-x|$ and $|z-x|$ are bounded below by $|x-p|$. Using these inequalities, we calculate
    \begin{align}\label{W_integral}
        \int_W\frac{\delta(x)^{2s}}{|y-x||z-x|}dx&\leq\int_{(p,p+M]}\frac{|x-p|^{2s-1}}{M}dx+\int_{[p+M,p+R)}|x-p|^{2s-2}dx.
    \end{align}
    When $s=1/2$, inequality (\ref{W_integral}) simplifies to
    \begin{align}\label{s=1/2_W_integral}
        \int_W\frac{\delta(x)^{2s}}{|y-x||z-x|}dx&\leq\int_{(p,p+M]}\frac{1}{M}dx+\int_{[p+M,p+R)}\frac{1}{|x-p|}dx\nonumber\\
        &=\frac{1}{M}\cdot M+\log|x-p|\Big|_{p+M}^{p+R}\nonumber\\
        &=1 + \log\left(\frac{R}{M}\right).
    \end{align}
    Whereas for $s<1/2$, we get
    \begin{align}\label{s<1/2_W_integral}
        \int_W\frac{\delta(x)^{2s}}{|y-x||z-x|}dx&\leq\int_{(p,p+M]}\frac{|x-p|^{2s-1}}{M}dx+\int_{[p+M,p+R)}|x-p|^{2s-2}dx \nonumber\\
        &\leq \frac{1}{M}\frac{1}{2s}(x-p)^{2s}\Big|_p^{p+M}+\frac{1}{2s-1}(x-p)^{2s-1}\Big|_{p+M}^{p+R}\nonumber\\
        &\leq C_sM^{2s-1},
    \end{align}
    where $C_s$ is a positive constant depending on $s$. So in this case, $Q_s(\cdot,\cdot)$ has the following upper bounds:
    \begin{align}\label{case4_Qbound_s<=1/2}
       Q_s(y,z)&\leq\begin{cases}
           \frac{C^2}{\delta(y)^s\delta(z)^s}\left(1+\log\left(\frac{R}{M}\right)\right), \qquad &s=1/2,\qquad y,z\in I\setminus\overline{U}\nonumber\\
           \frac{C^2}{\delta(y)^s\delta(z)^s}\left(1+M^{2s-1}\right), \qquad &s<1/2, \qquad y,z\in I\setminus\overline{U}
       \end{cases}\\
       &\leq\begin{cases}
            \frac{C^2}{\delta(y)^{1/2}\delta(z)^{1/2}}\left(1+\sum_{p\in\partial U}\left|\log\max\{|z-p|,|y-p|\}\right|\right), \qquad &s=1/2\\
            \frac{C^2}{\delta(y)^s\delta(z)^s}\left(1+\sum_{p\in\partial U}\max\{|y-p|,|z-p|\}^{2s-1}\right),\qquad &s<1/2.
        \end{cases}
    \end{align}
    On $A:=((I\setminus\overline{U}))\times (I\setminus\overline{U})))\setminus\left(\bigcup_{p\in\partial U}(B_R(p)\times B_R(p))\right)$, the terms in brackets are bounded.  Thus, we find that $Q_s(\cdot,\cdot)^q$ is integrable for $q<1/s$, just like in the previous cases.
    On $B_p=\left((I\setminus\overline{U})\cap B_R(p)\right)\times\left((I\setminus\overline{U})\cap B_R(p)\right)$ for each $p\in\partial U$, the integrability of $Q_s(\cdot,\cdot)^q$ is as follows: For $s<1/2$, set $u:=|y-p|$ and $v:=|z-p|$. The parameters $u$ and $v$ vary from $0$ to $R$ since $y$ and $z$ are both within a distance $R$ from $p$. W.l.o.g assume $v=\max\{u,v\}$. We have
    \begin{align}
        \int_{B_p}\frac{\max\{|y-p|,|y-p|\}^{(2s-1)q}}{u^{sq}v^{sq}}dydz&\leq C^{2q}\int_0^R\int_0^R\frac{\max\{u,v\}^{(2s-1)q}}{u^{sq}v^{sq}}dudv\nonumber\\
        &\leq 2C^{2q}\int_0^R u^{-sq}\int_u^R\frac{v^{(2s-1)q}}{v^{sq}}dvdu\nonumber\\
        &\leq 2C^{2q}\int_0^R u^{-sq}\int_u^R v^{q(s-1)}dvdu\nonumber\\
        &\leq 2C^{2q}\int_0^R u^{-sq}\left(\frac{1}{q(s-1)+1} v^{q(s-1)+1}\Big|_u^R\right)du\nonumber\\
        &=\frac{2C^{2q}}{q(s-1)+1}\int_0^R u^{-sq}\left(D_{q,s}-u^{q(s-1)+1}\right)du\nonumber\\
        &=\frac{2C^{2q}}{q(s-1)+1}\int_0^R D_{q,s}u^{-sq}-u^{1-q}du\nonumber\\
        &<\infty,
    \end{align}
    since the first term in the last integral is integrable for $q<1/s$ and $s<1/2$, and the second term is integrable for $1-q>-1$ or $q<2$. $D_{q,s}$ is a sufficiently large positive constant depending on $q$ and $s$. Thus, for $s<1/2$, $Q_s(\cdot,\cdot)^q$ is integrable for $q<2$. For $s=1/2$, we can bound the log term by
    \begin{align*}
        |\log\max\{|y-p|,|z-p|\}|\leq C_\epsilon\max\{|y-p|,|z-p|\}^{-\epsilon},
    \end{align*}
    where $\epsilon>0$ is arbitrary and $C_\epsilon$ is a positive constant depending on $\epsilon$. Arguing as before, we find that $Q_s(\cdot,\cdot)^q$ is integrable when $-\epsilon q-\frac{q}{2}+1>-1$ or $q<\frac{2}{\epsilon+\frac{1}{2}}$. However, since $\epsilon$ is arbitrary, we conclude that integrabililty holds when $q<4$.
    
    Substituting the upper bounds for $Q_s(\cdot,\cdot)^q$ into inequality (\ref{L2_inquality_w}) yields an $L^p$ bound for $\|w\|^2_{L^2(U)}$. Thus we have
    \begin{align}\label{LP_bound_w} 
        \|w\|_{L^2(U)}\leq \tilde{C}\|u\|_{L^p(\R\setminus\overline{U})}=\tilde{C}\|w\|_{L^p(\R\setminus\overline{U})}
    \end{align}
    provided that $p>\frac{1}{1-s}$ for $s>1/2$, and $p>2$ for $s\leq 1/2$. $\tilde{C}$ is a large positive constant depending on $s$ and $U$.  
\end{proof}

$L^2$ convergence of the finite potential well eigenfunctions $u_{i,\delta}$ to the infinite potential well eigenfunctions $u_i$ follows from a compactness and energy minimisation argument within a proof by induction.
The energy functional associated with the operator in the infinite potential well eigenvalue problem (\ref{evalue_problem_infinite_well}) is
\begin{align}\label{energy_infinite_well}
    \mathcal{E}[u]=\frac{c_s}{2}[u]^2_{H^s(\R)}+\int_{U}V(x)u(x)^2dx
\end{align}
for any function $u$ in $H^s_0(U)$. Recall the semi-norm $[\cdot]_{H^s}$ is given by equation (\ref{def_Gagliardo_seminorm}) in Section \ref{subs_func_an} by taking $n=1$, i.e.,
\begin{align}\label{seminorm}
    [u]^2_{H^s(\R)}=\int_{\R}\int_{\R}\frac{|u(x)-u(y))|^2}{|x-y|^{1+2s}}dxdy.
\end{align}
Similarly, the energy functional associated with the operator in the finite potential well eigenvalue problem (\ref{evalue_problem_finite_well}) is
\begin{align}\label{energy_finite_well}
    \mathcal{E}_{\delta}[u]=\mathcal{E}[u]+\frac{1}{\delta}\int_{I\setminus U} u(x)^2dx=\mathcal{E}[u]+\frac{1}{\delta}\|u\|^2_{L^2(\R\setminus U)}
\end{align}
for any function $u$ in $H^s_0(I)$ and $\delta>0$. Clearly $\mathcal{E}[u]=\mathcal{E}_{\delta}[u]$ if $u$ belongs to $H^s_0(U)$.

Using equations \eqref{energy_infinite_well} and
\eqref{energy_finite_well}, the variational characterisations of
$\lambda_i$ and $\lambda_{i,\delta}$ are, for $\delta>0$ and
$i=1,2,\ldots$,
\begin{align}\label{var_char_lambda_i}
    \lambda_i
    =
    \inf_{\substack{S\subset H_0^s(U)\\ \dim S=i}}
    \;
    \sup_{\substack{u\in S\\ \|u\|_{L^2(U)}=1}}
    \mathcal E[u]
\end{align}
and
\begin{align}\label{var_char_lambda_idelta}
    \lambda_{i,\delta}
    =
    \inf_{\substack{S\subset H_0^s(I)\\ \dim S=i}}
    \;
    \sup_{\substack{u\in S\\ \|u\|_{L^2(I)}=1}}
    \mathcal E_\delta[u].
\end{align}

\begin{lemma}\label{lem_eigenvalue_inequality}
Fix $\delta>0$, and let $(u_i,\lambda_i)$ and
$(u_{i,\delta},\lambda_{i,\delta})$ be the $i$th eigenpairs satisfying
eigenvalue problems \eqref{evalue_problem_infinite_well} and
\eqref{evalue_problem_finite_well}, respectively. Then
\begin{align}\label{eigenvalue_inequality}
    \lambda_{i,\delta}\leq\lambda_i.
\end{align}
\end{lemma}

\begin{proof}
Since $H_0^s(U)\subset H_0^s(I)$, we may restrict the infimum in
\eqref{var_char_lambda_idelta} to $i$-dimensional subspaces of
$H_0^s(U)$. Thus
\[
    \lambda_{i,\delta}
    \leq
    \inf_{\substack{S\subset H_0^s(U)\\ \dim S=i}}
    \;
    \sup_{\substack{u\in S\\ \|u\|_{L^2(I)}=1}}
    \mathcal E_\delta[u].
\]
For $u\in H_0^s(U)$, we have $u=0$ on $I\setminus U$, and hence
\[
    \|u\|_{L^2(I)}
    =
    \|u\|_{L^2(U)}
    \qquad\text{and}\qquad
    \mathcal E_\delta[u]
    =
    \mathcal E[u].
\]
Therefore,
\[
    \lambda_{i,\delta}
    \leq
    \inf_{\substack{S\subset H_0^s(U)\\ \dim S=i}}
    \;
    \sup_{\substack{u\in S\\ \|u\|_{L^2(U)}=1}}
    \mathcal E[u]
    =
    \lambda_i.
\]
\end{proof}
\noindent A simple corollary is that $u_{_{i,\delta}}$ is small in $L^2$ on $\R\setminus\overline{U}$.
\begin{corollary}\label{small_L2_bound_w_Ucomp}
Fix $0<\delta\ll 1$, and let $u_{i,\delta}$ be an eigenfunction of the
$i$th eigenvalue $\lambda_{i,\delta}$ of eigenvalue problem
\eqref{evalue_problem_finite_well}. Then
\begin{align}
    \|u_{i,\delta}\|_{L^2(\R\setminus\overline{U})}
    \leq C\delta^{1/2},
\end{align}
where $C>0$ is independent of $\delta$.
\end{corollary}

\begin{proof}
By Lemma~\ref{lem_eigenvalue_inequality},
\[
    \lambda_{i,\delta}\leq\lambda_i.
\]
Since the eigenfunctions $u_{i,\delta}$ are normalised in $L^2(I)$,
we have
\[
    \lambda_{i,\delta}
    =
    \mathcal E_\delta[u_{i,\delta}],
\]
which implies that
\begin{align*}
    \lambda_i
    &\geq \lambda_{i,\delta}\\
    &=
    \frac{c_s}{2}
    [u_{i,\delta}]_{H^s(\R)}^2
    +
    \int_U V(x)u_{i,\delta}(x)^2\,dx
    +
    \frac{1}{\delta}
    \|u_{i,\delta}\|_{L^2(I\setminus U)}^2\\
    &\geq
    \frac{1}{\delta}
    \|u_{i,\delta}\|_{L^2(I\setminus U)}^2.
\end{align*}
The last inequality follows from the non-negativity of the fractional Dirichlet energy and the potential $V$.

It follows that
\[
    \|u_{i,\delta}\|_{L^2(I\setminus U)}^2
    \leq
    \lambda_i\delta,
\]
and hence
\[
    \|u_{i,\delta}\|_{L^2(I\setminus U)}
    \leq
    \sqrt{\lambda_i}\,\delta^{1/2}.
\]
Since $u_{i,\delta}=0$ on $\R\setminus I$ and $\partial U$ has
measure zero,
\[
    \|u_{i,\delta}\|_{L^2(\R\setminus\overline U)}
    =
    \|u_{i,\delta}\|_{L^2(I\setminus U)}
    \leq
    C\delta^{1/2},
\]
where we may take $C=\sqrt{\lambda_i}$.
\end{proof}

\noindent A natural second corollary is the following.
\begin{corollary}\label{asymptotic_norm_orthogonality_U}
Let $u_{i,\delta}$ be the $L^2(I)$-orthonormal eigenfunctions of
eigenvalue problem \eqref{evalue_problem_finite_well}. Then
\begin{align}\label{L2_norm_on_U}
    \sqrt{1-C^2\delta}
    \leq
    \|u_{i,\delta}\|_{L^2(U)}
    \leq1.
\end{align}
In particular,
\[
    \|u_{i,\delta}\|_{L^2(U)}
    \longrightarrow1
    \qquad\text{as }\delta\to0^+.
\]
Moreover, for $i\neq j$,
\begin{align}\label{asymptotic_orthogonality_U}
    \left|
    \langle u_{i,\delta},u_{j,\delta}\rangle_{L^2(U)}
    \right|
    \leq C_{i,j}\delta.
\end{align}
\end{corollary}

\begin{proof}
Since $\|u_{i,\delta}\|_{L^2(I)}=1$,
\[
    \|u_{i,\delta}\|_{L^2(U)}^2
    =
    1-\|u_{i,\delta}\|_{L^2(I\setminus U)}^2.
\]
Noting that $\|u_{i,\delta}\|_{L^2(\R\setminus\overline{U})}=
    \|u_{i,\delta}\|_{L^2(I\setminus U)}$, Corollary~\ref{small_L2_bound_w_Ucomp} implies that
\[
    1-C^2\delta
    \leq
    \|u_{i,\delta}\|_{L^2(U)}^2
    \leq1,
\]
which proves \eqref{L2_norm_on_U}. In particular,
\[
    \|u_{i,\delta}\|_{L^2(U)}
    \longrightarrow1
    \qquad\text{as }\delta\to0^+.
\]

For $i\neq j$, orthogonality in $L^2(I)$ yields
\[
    \langle u_{i,\delta},u_{j,\delta}\rangle_{L^2(U)}
    =
    -
    \langle u_{i,\delta},u_{j,\delta}\rangle_{L^2(I\setminus U)}.
\]
Applying Corollary~\ref{small_L2_bound_w_Ucomp} again gives
\begin{align*}
    \left|
    \langle u_{i,\delta},u_{j,\delta}\rangle_{L^2(U)}
    \right|
    &\leq
    \|u_{i,\delta}\|_{L^2(I\setminus U)}
    \|u_{j,\delta}\|_{L^2(I\setminus U)}\\
    &\leq C_{i,j}\delta.
\end{align*}
\end{proof}

We require one final energy estimate and an asymptotic orthogonality
property to prove our main result.
{\begin{lemma}\label{lem_refined_evalue_inequality}
For each fixed $i\geq1$, there exist constants $C>0$ and
$\eta,\rho>0$, independent of $\delta$, such that, for all
sufficiently small $\delta>0$,
\begin{align}\label{final_evalue_inequality}
\mathcal E[\hat v_{i,\delta}]
    \leq
    \frac{\lambda_i+C\delta^\rho}
         {(1-C\delta^\eta)^2},
\end{align}
where
    $\hat v_{i,\delta}
    :=
    \frac{v_{i,\delta}}
         {\|v_{i,\delta}\|_{L^2(U)}}$.
Moreover, for $i\neq j$,
\begin{align}\label{asymptotic_orthogonality_vhat}
    \left|
    \langle
    \hat v_{i,\delta},
    \hat v_{j,\delta}
    \rangle_{L^2(U)}
    \right|
    \longrightarrow0
    \qquad\text{as }\delta\to0^+.
\end{align}
\end{lemma}
\begin{proof}
    Using Lemma \ref{lem_eigenvalue_inequality}, the energy equations (\ref{energy_infinite_well}) and (\ref{energy_finite_well}), and recalling that $u_{i,\delta}=w_{i,\delta}+v_{i,\delta}$ with
    $\mathcal E_{\mathrm{res}}
    [v_{i,\delta},w_{i,\delta}]=0$,
we calculate: for $0<\delta\ll 1$ and $i=1,2,\ldots,$
\begin{align*}
    \lambda_i
    &\geq \lambda_{i,\delta}\\
    &=
    \frac{c_s}{2}[u_{i,\delta}]^2_{H^s(\R)}
    +\int_U V(x)u_{i,\delta}(x)^2\,dx
    +\frac{1}{\delta}
     \int_{I\setminus U}u_{i,\delta}(x)^2\,dx\\
    &=
    \frac{c_s}{2}[w_{i,\delta}+v_{i,\delta}]^2_{H^s(\R)}
    +\int_U V(x)
       \left(w_{i,\delta}+v_{i,\delta}\right)^2\,dx
    +\frac{1}{\delta}
     \int_{I\setminus U}u_{i,\delta}(x)^2\,dx\\
    &=
    \mathcal E[v_{i,\delta}]
    +\mathcal E_\delta[w_{i,\delta}]
    +2\int_U
      V(x)w_{i,\delta}(x)v_{i,\delta}(x)\,dx\\
    &\geq
    \mathcal E[v_{i,\delta}]
    -2\|V\|_{L^\infty(U)}
      \|w_{i,\delta}\|_{L^2(U)}
      \|v_{i,\delta}\|_{L^2(U)}.
\end{align*} 
Thus,
\begin{align}\label{evalue_inequality_induction}
    \lambda_i
    \geq\lambda_{i,\delta}
    \geq
    \mathcal E[v_{i,\delta}]
    -2\|V\|_{L^\infty(U)}
      \|w_{i,\delta}\|_{L^2(U)}
      \|v_{i,\delta}\|_{L^2(U)}.
\end{align}


From the first and second lines of the above calculation, it can also be deduced that
\[
    \frac{c_s}{2}
    [u_{i,\delta}]_{H^s(\R)}^2
    \leq
    \lambda_{i,\delta}
    \leq
    \lambda_i.
\]
Recalling
\[
    \|u\|_{H^s(\R)}^2
    =
    \|u\|_{L^2(\R)}^2
    +
    \frac{c_s}{2}[u]_{H^s(\R)}^2,
\]
and using $\|u_{i,\delta}\|_{L^2(\R)}
=\|u_{i,\delta}\|_{L^2(I)}=1$, we obtain
\[
    \|u_{i,\delta}\|_{H^s(\R)}^2
    \leq 1+\lambda_i.
\]
Consequently,
\[
    \|u_{i,\delta}\|_{H^s(\R)}
    \leq C_i
\]
uniformly for sufficiently small $\delta>0$.


We now show that $w_{i,\delta}$ can be made arbitrarily small in $L^2(U)$. When $s>1/2$, we apply Theorem 8.2 in \cite{NPV} with
$\Omega=\R$, $n=1$, and $p=2$. 
Since $2s>1$, Theorem 8.2 and the preceding uniform
$H^s(\R)$ bound give
\[
    \|u_{i,\delta}\|_{C^{0,s-1/2}(\R)}
    \leq C_i
\]
uniformly for sufficiently small $\delta>0$.
In particular,
\[
    \|u_{i,\delta}\|_{L^\infty(\R\setminus\overline U)}
    \leq
    \|u_{i,\delta}\|_{L^\infty(\R)}
    \leq
    \|u_{i,\delta}\|_{C^{0,s-1/2}(\R)}
    \leq C_i.
\]
It follows from an interpolation inequality, Corollary \ref{small_L2_bound_w_Ucomp}, and the above
$L^\infty$ bound that,
for every fixed finite $p\geq 2$,

\begin{align}\label{1_Lq_bound}
\|u_{i,\delta}\|_{L^p(\R\setminus\overline U)}
&\leq
\|u_{i,\delta}\|_{L^2(\R\setminus\overline U)}^{2/p}
\|u_{i,\delta}\|_{L^\infty(\R\setminus\overline U)}^{1-2/p}\leq C\delta^{1/p}.
\end{align}

If $s<1/2$, the fractional Sobolev inequality in
\cite[Th.~6.7]{NPV}, applied on $\R$ with $n=1$ and $p=2$, gives
\[
    \|u_{i,\delta}\|_{L^q(\R)}
    \leq
    C\|u_{i,\delta}\|_{H^s(\R)}
    \leq C_i
\]
for
\[
    2\leq q\leq\frac{2}{1-2s}.
\]
Therefore,
\[
    \|u_{i,\delta}\|_{L^q(\R\setminus\overline U)}
    \leq C_i.
\]
When $s=1/2$, the critical fractional Sobolev inequality
\cite[Th.~6.10]{NPV} gives the same estimate for every finite
$q\geq2$.

Now interpolate this uniform $L^q$ bound and apply the estimate from Corollary~\ref{small_L2_bound_w_Ucomp}. For $2<p<q$, let
$\theta\in(0,1)$ be determined by
\[
    \frac1p
    =
    \frac{\theta}{2}
    +
    \frac{1-\theta}{q}.
\]
Then
\begin{align*}
    \|u_{i,\delta}\|_{L^p(\R\setminus\overline U)}
    &\leq
    \|u_{i,\delta}\|_{L^2(\R\setminus\overline U)}^\theta
    \|u_{i,\delta}\|_{L^q(\R\setminus\overline U)}^{1-\theta}\leq
    C\delta^{\theta/2}.
\end{align*}
If $s<1/2$, choose
\[
    q=\frac{2}{1-2s}.
\]
Then
\[
    \theta
    =
    \frac{2-p(1-2s)}{2sp},
\]
and 
\begin{align}\label{2_Lq_bound}
    \|u_{i,\delta}\|_{L^p(\R\setminus\overline U)}
    \leq
    C\delta^{\frac{2-p(1-2s)}{4sp}},
    \qquad
    2<p<\frac{2}{1-2s}.
\end{align}
If $s=1/2$, choose any finite $q>p>2$. Then
\[
    \theta
    =
    \frac{2(q-p)}{p(q-2)},
\]
and 
\begin{align}\label{3_Lq_bound}
    \|u_{i,\delta}\|_{L^p(\R\setminus\overline U)}
    \leq
    C\delta^{\frac{q-p}{p(q-2)}}.
\end{align}


Therefore, since $w_{i,\delta}=u_{i,\delta}$ on
$\R\setminus U$, Lemma~\ref{L2_estimate_w} and the foregoing
$L^p$ estimates give
\begin{align}\label{small_L2_bound_w}
    \|w_{i,\delta}\|_{L^2(U)}
    &\leq
    \widetilde C
    \|w_{i,\delta}\|_{L^p(\R\setminus\overline U)}=
    \widetilde C
    \|u_{i,\delta}\|_{L^p(\R\setminus\overline U)}
    \leq
    C'\delta^\rho,
\end{align}
where, in each of the three cases above, $p$ is chosen to satisfy
the corresponding condition in Lemma~\ref{L2_estimate_w}, and
$\rho>0$ is the resulting exponent of $\delta$.

Combining the reverse triangle inequality with inequalities (\ref{small_L2_bound_w}) and
(\ref{L2_norm_on_U}) yields
\begin{align*}
    \|v_{i,\delta}\|_{L^2(U)}
    &=\|u_{i,\delta}-w_{i,\delta}\|_{L^2(U)}\\
    &\geq
    \|u_{i,\delta}\|_{L^2(U)}
    -\|w_{i,\delta}\|_{L^2(U)}\\
    &\geq\sqrt{1-C^2\delta}-C'\delta^\rho\\
    &\geq 1-C^2\delta-C'\delta^\rho,\qquad 0<\delta\ll 1
\end{align*}
Let $\eta:=\min\{\rho,1\}$. After relabelling the constants, we have
\begin{align}\label{lower_bound_v}
    \|v_{i,\delta}\|_{L^2(U)}
    \geq1-C\delta^\eta.
\end{align}
On the other hand, recalling that $\|u_{i,\delta}\|_{L^2(I)}=1$,
\[
    \|v_{i,\delta}\|_{L^2(U)}
    \leq
    \|u_{i,\delta}\|_{L^2(U)}
    +\|w_{i,\delta}\|_{L^2(U)}
    \leq1+C\delta^\rho.
\]
Consequently,
\[
    \|v_{i,\delta}\|_{L^2(U)}\longrightarrow1\qquad\mbox{as}\qquad\delta\to 0^+.
\]
In particular, $v_{i,\delta}\neq0$ for sufficiently small $\delta$,
so
\[
    \hat v_{i,\delta}
    =
    \frac{v_{i,\delta}}
         {\|v_{i,\delta}\|_{L^2(U)}}
\]
is well defined. Using (\ref{evalue_inequality_induction}),
(\ref{small_L2_bound_w}), and the preceding upper bound for
$\|v_{i,\delta}\|_{L^2(U)}$, we obtain
\[
    \lambda_i
    \geq\lambda_{i,\delta}
    \geq
    \mathcal E[\hat v_{i,\delta}]
    \|v_{i,\delta}\|_{L^2(U)}^2
    -C\delta^\rho.
\]
Since $\mathcal E[\hat v_{i,\delta}]\geq0$, inequality
(\ref{lower_bound_v}) gives
\begin{align*}
    \lambda_i
    \geq\lambda_{i,\delta}
    \geq
    \mathcal E[\hat v_{i,\delta}]
    (1-C\delta^\eta)^2
    -C\delta^\rho.
\end{align*}
Rearranging this inequality for $\mathcal{E}[\hat{v}_{i,\delta}]$ leads to inequality~(\ref{final_evalue_inequality}).

Finally we show that the normalised projections are asymptotically orthogonal. For $i\neq j$, equation
(\ref{asymptotic_orthogonality_U}), the identity
$v_{m,\delta}=u_{m,\delta}-w_{m,\delta}$, and
(\ref{small_L2_bound_w}) give
\begin{align*}
\left|
\langle v_{i,\delta},v_{j,\delta}\rangle_{L^2(U)}
\right|
&\leq
\left|
\langle u_{i,\delta},u_{j,\delta}\rangle_{L^2(U)}
\right|\\
&\quad+
\|w_{i,\delta}\|_{L^2(U)}
\|u_{j,\delta}\|_{L^2(U)}
+
\|u_{i,\delta}\|_{L^2(U)}
\|w_{j,\delta}\|_{L^2(U)}\\
&\quad+
\|w_{i,\delta}\|_{L^2(U)}
\|w_{j,\delta}\|_{L^2(U)}\\
&\leq
C\bigl(\delta+\delta^\rho+\delta^{2\rho}\bigr).
\end{align*}
Together with (\ref{lower_bound_v}), this yields
\[
    \left|
    \langle
    \hat v_{i,\delta},
    \hat v_{j,\delta}
    \rangle_{L^2(U)}
    \right|
    \longrightarrow0
    \qquad
    \text{as }\delta\to0^+.
\]
\end{proof}

\begin{theorem}\label{L2_convergence_efunctions} There exists a sequence of deltas $\delta$ approaching zero such that $u_{i,\delta}$ converges in $L^2$ to $u_i$ for each $i$, where $u_i$ is an eigenfunction for eigenvalue $\lambda_i$, and the sequence of eigenfunctions $\{u_i\}_{i\geq 1}$ is an orthonormal basis for $L^2(U)$.
\end{theorem}
\begin{proof} Both $u_i$ and $u_{i,\delta}$ are zero off the interval $I=(-1,1)$, so the claim is trivial in this case. On the set $I\setminus U$, $u_i$ is still zero by definition, while $\|u_{i,\delta}\|_{L^2(I\setminus U)}\to 0$ as $\delta\to 0^+$ by
Corollary \ref{small_L2_bound_w_Ucomp}.
This implies that
\[
    \|u_i-u_{i,\delta}\|_{L^2(I\setminus U)}
    =
    \|u_{i,\delta}\|_{L^2(I\setminus U)}
    \to0
\]
as $\delta\to0^+$.

The only case that remains to be proven is $L^2$ convergence of $u_{i,\delta}$ to $u_i$ on $U$ as $\delta\to0^+$. We prove this by induction, constructing a nested sequence of subsequences and then taking a diagonal sequence along which the required convergence holds for every $i$. 

\textbf{Compactness argument for fixed $i$.}
Fix $i\geq1$. For sufficiently small $\delta>0$,
Lemma~\ref{lem_refined_evalue_inequality} gives
\[
    \mathcal E[\hat v_{i,\delta}]
    \leq
    \frac{\lambda_i+C\delta^\rho}
         {(1-C\delta^\eta)^2}.
\]
The right-hand side is uniformly bounded for sufficiently
small $\delta>0$, so there exists $M_i>0$, independent of $\delta$,
such that
\[
    \mathcal E[\hat v_{i,\delta}]\leq M_i.
\]
As $V\geq0$ and
\[
    \mathcal E[\hat v_{i,\delta}]
    =
    \mathcal E_{\mathrm{res}}
    [\hat v_{i,\delta},\hat v_{i,\delta}]
    +
    \int_U V(x)\hat v_{i,\delta}(x)^2\,dx,
\]
we conclude
\[
    \|\hat v_{i,\delta}\|_{H_0^s(U)}
    =
    \mathcal E_{\mathrm{res}}
    [\hat v_{i,\delta},\hat v_{i,\delta}]^{1/2}
    \leq \sqrt{M_i}.
\]
Thus $\{\hat v_{i,\delta}\}$ is bounded in $H_0^s(U)$.

Choose any sequence $\delta_n\to0$. Since $H_0^s(U)$ is a Hilbert
space, it is reflexive \cite[Th.~4.6-6]{Krey}. By
\cite[Th.~3.18]{Brezis}, after passing to a subsequence, and
retaining the notation $\delta_n$, there exists $u\in H_0^s(U)$
such that
\[
    \hat v_{i,\delta_n}
    \rightharpoonup u
    \qquad\text{in }H_0^s(U).
\]
In terms of the inner product used on $H_0^s(U)$, this means
that
\[
    \mathcal E_{\mathrm{res}}
    [\hat v_{i,\delta_n},\varphi]
    \longrightarrow
    \mathcal E_{\mathrm{res}}[u,\varphi]
    \qquad
    \text{for every }\varphi\in H_0^s(U).
\]

Next we obtain convergence in $L^2(U)$. As $U$ is a finite union
of bounded intervals, it is a bounded Lipschitz open set and hence
an extension domain for $H^s(U)$ by \cite[Th.~5.4]{NPV}, using the
identification $H^s=W^{s,2}$. Furthermore, since
\[
    \|\hat v_{i,\delta_n}\|_{L^2(U)}=1
    \qquad\text{for every }n,
\]
the sequence $\{\hat v_{i,\delta_n}\}_{n\geq1}$ is uniformly
bounded in $L^2(U)$. The preceding $H_0^s(U)$ bound also implies
that
\[
    \sup_{n\geq1}
    \iint_{U\times U}
    \frac{
        |\hat v_{i,\delta_n}(x)
        -\hat v_{i,\delta_n}(y)|^2
    }{
        |x-y|^{1+2s}
    }
    \,dx\,dy
    <\infty.
\]
It follows from Theorem~7.1 in \cite{NPV}, with $p=q=2$, that
$\{\hat v_{i,\delta_n}\}_{n\geq1}$ is precompact in $L^2(U)$.
Passing to a further subsequence if necessary, and retaining the
notation $\delta_n$, there exists $\widetilde u\in L^2(U)$ such
that
\[
    \hat v_{i,\delta_n}
    \longrightarrow \widetilde u
    \qquad\text{in }L^2(U).
\]

By Proposition~\ref{prop_equiv_norms}, $H_0^s(U)$ is continuously
embedded in $L^2(U)$, so weak convergence in $H_0^s(U)$ implies
weak convergence in $L^2(U)$; thus
\[
    \hat v_{i,\delta_n}
    \rightharpoonup u
    \qquad\text{in }L^2(U).
\]
On the other hand, $L^2(U)$ convergence to $\widetilde u$ also
implies weak $L^2(U)$ convergence to $\widetilde u$. By uniqueness
of weak limits \cite[Lem.~4.8-3]{Krey},
\[
    \widetilde u=u.
\]
So
\[
    \hat v_{i,\delta_n}
    \rightharpoonup u
    \qquad\text{weakly in }H_0^s(U),
\]
and
\[
    \hat v_{i,\delta_n}
    \longrightarrow u
    \qquad\text{strongly in }L^2(U).
\]

Finally, again using the fact that
$\|\hat v_{i,\delta_n}\|_{L^2(U)}=1$ for every $n$, the reverse
triangle inequality gives
\[
    \left|
        \|\hat v_{i,\delta_n}\|_{L^2(U)}
        -
        \|u\|_{L^2(U)}
    \right|
    \leq
    \|\hat v_{i,\delta_n}-u\|_{L^2(U)}
    \longrightarrow0,
\]
and hence
\[
    \|u\|_{L^2(U)}=1.
\]

The $H_0^s(U)$ norm is weakly lower semicontinuous
\cite[Sec.~4.8, Prob.~4]{Krey}, so
\[
    \mathcal E_{\mathrm{res}}[u,u]
    \leq
    \liminf_{n\to\infty}
    \mathcal E_{\mathrm{res}}
    [\hat v_{i,\delta_n},\hat v_{i,\delta_n}].
\]
Moreover, since $\hat v_{i,\delta_n}\to u$ in $L^2(U)$ and
$V\in L^\infty(U)$,
\[
    \int_U V(x)\hat v_{i,\delta_n}(x)^2\,dx
    \longrightarrow
    \int_U V(x)u(x)^2\,dx.
\]
Therefore,
\begin{align*}
    \mathcal E[u]
    &=
    \mathcal E_{\mathrm{res}}[u,u]
    +
    \int_U V(x)u(x)^2\,dx\\
    &\leq
    \liminf_{n\to\infty}
    \mathcal E_{\mathrm{res}}
    [\hat v_{i,\delta_n},\hat v_{i,\delta_n}]
    +
    \lim_{n\to\infty}
    \int_U V(x)\hat v_{i,\delta_n}(x)^2\,dx\\
    &=
    \liminf_{n\to\infty}
    \left(
        \mathcal E_{\mathrm{res}}
        [\hat v_{i,\delta_n},\hat v_{i,\delta_n}]
        +
        \int_U V(x)\hat v_{i,\delta_n}(x)^2\,dx
    \right)\\
    &=
    \liminf_{n\to\infty}
    \mathcal E[\hat v_{i,\delta_n}].
\end{align*}
Lemma~\ref{lem_refined_evalue_inequality} also gives
\[
    \mathcal E[\hat v_{i,\delta_n}]
    \leq
    \frac{\lambda_i+C\delta_n^\rho}
         {(1-C\delta_n^\eta)^2}
\]
for all sufficiently large $n$. Taking the $\liminf$ on both sides
and using $\delta_n\to0$ gives
\[
    \liminf_{n\to\infty}
    \mathcal E[\hat v_{i,\delta_n}]
    \leq
    \liminf_{n\to\infty}
    \frac{\lambda_i+C\delta_n^\rho}
         {(1-C\delta_n^\eta)^2}
    =
    \lambda_i.
\]
Consequently,
\begin{align}\label{energy_bound_subseq_limit}
    \mathcal E[u]\leq\lambda_i.
\end{align}
We now use this compactness argument to construct the eigenfunctions
and the required nested subsequences.

\textbf{Base case.}
Choose any sequence $\delta_n\to0$. Applying the preceding argument
with $i=1$, after passing to a subsequence, which we denote by
$\delta_n^{(1)}$, there exists $u_1\in H_0^s(U)$ such that
\[
    \hat v_{1,\delta_n^{(1)}}
    \rightharpoonup u_1
    \qquad\text{in }H_0^s(U),
\]
and
\[
    \hat v_{1,\delta_n^{(1)}}
    \longrightarrow u_1
    \qquad\text{in }L^2(U),
\]
with
\[
    \|u_1\|_{L^2(U)}=1
    \qquad\text{and}\qquad
    \mathcal E[u_1]\leq\lambda_1.
\]
On the other hand, the variational characterisation of $\lambda_1$
gives
\[
    \mathcal E[u_1]\geq\lambda_1.
\]
It follows that
\[
    \mathcal E[u_1]=\lambda_1,
\]
so $u_1$ is an eigenfunction corresponding to $\lambda_1$.

\textbf{Induction step.}
Suppose inductively that we have constructed nested subsequences
\[
    \{\delta_n^{(1)}\}
    \supset
    \{\delta_n^{(2)}\}
    \supset\cdots\supset
    \{\delta_n^{(k)}\},
    \qquad
    \delta_n^{(k)}\to0,\qquad n\to\infty,
\]
and orthonormal eigenfunctions $u_1,\ldots,u_k$, corresponding
respectively to $\lambda_1,\ldots,\lambda_k$, such that
\[
    \hat v_{j,\delta_n^{(k)}}
    \longrightarrow u_j
    \qquad\text{in }L^2(U),
    \qquad j=1,\ldots,k.
\]

As $\delta_n^{(k)}\to0$, we may apply the preceding compactness
argument with $i=k+1$ to the sequence $\{\delta_n^{(k)}\}$.
After passing to a further subsequence
\[
    \{\delta_n^{(k+1)}\}
    \subset
    \{\delta_n^{(k)}\},
\]
there exists $u_{k+1}\in H_0^s(U)$ such that
\[
    \hat v_{k+1,\delta_n^{(k+1)}}
    \rightharpoonup u_{k+1}
    \quad\text{in }H_0^s(U),
\]
and
\[
    \hat v_{k+1,\delta_n^{(k+1)}}
    \longrightarrow u_{k+1}
    \quad\text{in }L^2(U),
\]
with
\[
    \|u_{k+1}\|_{L^2(U)}=1
    \qquad\text{and}\qquad
    \mathcal E[u_{k+1}]
    \leq\lambda_{k+1}.
\]

Since $\{\delta_n^{(k+1)}\}\subset\{\delta_n^{(k)}\}$, the convergence already established for
$j=1,\ldots,k$ is preserved along this further subsequence. As a result, for each $j=1,\ldots,k$,
\begin{align*}
&\left|
    \langle u_{k+1},u_j\rangle_{L^2(U)}
    -
    \left\langle
        \hat v_{k+1,\delta_n^{(k+1)}},
        \hat v_{j,\delta_n^{(k+1)}}
    \right\rangle_{L^2(U)}
\right|\\
&\qquad\leq
    \left\|
        u_{k+1}
        -
        \hat v_{k+1,\delta_n^{(k+1)}}
    \right\|_{L^2(U)}
    \|u_j\|_{L^2(U)}
\\
&\qquad\quad+
    \left\|
        \hat v_{k+1,\delta_n^{(k+1)}}
    \right\|_{L^2(U)}
    \left\|
        u_j
        -
        \hat v_{j,\delta_n^{(k+1)}}
    \right\|_{L^2(U)}
    \longrightarrow0.
\end{align*}
Hence,
\[
    \langle u_{k+1},u_j\rangle_{L^2(U)}
    =
    \lim_{n\to\infty}
    \left\langle
        \hat v_{k+1,\delta_n^{(k+1)}},
        \hat v_{j,\delta_n^{(k+1)}}
    \right\rangle_{L^2(U)}
    =0,
\]
where the final equality follows from the asymptotic orthogonality
in Lemma~\ref{lem_refined_evalue_inequality}, as
$\delta_n^{(k+1)}\to0$. It follows that
\[
    u_{k+1}\perp u_j,
    \qquad j=1,\ldots,k,
\]
so $u_{k+1}$ is a unit vector orthogonal to
$\operatorname{span}\{u_1,\ldots,u_k\}$.

On the other hand, the variational characterisation
\[
    \lambda_{k+1}
    =
    \inf_{\substack{
        u\in H_0^s(U),\ \|u\|_{L^2(U)}=1\\
        u\perp\operatorname{span}\{u_1,\ldots,u_k\}
    }}
    \mathcal E[u]
\]
implies, as $u_{k+1}$ is a unit vector orthogonal to
$u_1,\ldots,u_k$, that
\[
    \mathcal E[u_{k+1}]\geq\lambda_{k+1}.
\]
Together with the earlier upper bound,
\[
    \mathcal E[u_{k+1}]
    \leq\lambda_{k+1},
\]
this gives
\[
    \mathcal E[u_{k+1}]
    =
    \lambda_{k+1},
\]
and $u_{k+1}$ is an eigenfunction corresponding to
$\lambda_{k+1}$. This argument also applies when eigenvalues have
multiplicity greater than one, as the eigenvalues are listed
according to multiplicity.

Repeating this construction produces a nested family of subsequences.
Taking the diagonal sequence
\[
    \delta_n:=\delta_n^{(n)},
\]
we have $\delta_n\to0$. For every fixed $i$, the tail of this
diagonal sequence from $n\geq i$ is contained in the nested
subsequences on which convergence of the $i$th projected
eigenfunction has already been demonstrated. Therefore,
\[
    \hat v_{i,\delta_n}
    \longrightarrow u_i
    \qquad\text{in }L^2(U)
\]
for every fixed $i$.

It remains to pass from $\hat v_{i,\delta_n}$ to
$u_{i,\delta_n}$. Recalling that
\[
    u_{i,\delta}=v_{i,\delta}+w_{i,\delta}
\]
and
\[
    v_{i,\delta}
    =
    \|v_{i,\delta}\|_{L^2(U)}\hat v_{i,\delta},
\]
we have
\begin{align*}
    \|u_{i,\delta}-\hat v_{i,\delta}\|_{L^2(U)}
    &\leq
    \|w_{i,\delta}\|_{L^2(U)}
    +
    \|v_{i,\delta}-\hat v_{i,\delta}\|_{L^2(U)}\\
    &=
    \|w_{i,\delta}\|_{L^2(U)}
    +
    \left|
        \|v_{i,\delta}\|_{L^2(U)}-1
    \right|.
\end{align*}
By \eqref{small_L2_bound_w},
\[
    \|w_{i,\delta}\|_{L^2(U)}\longrightarrow0.
\]
Also,
\[
    \left|
        \|v_{i,\delta}\|_{L^2(U)}
        -
        \|u_{i,\delta}\|_{L^2(U)}
    \right|
    \leq
    \|w_{i,\delta}\|_{L^2(U)},
\]
while Corollary~\ref{asymptotic_norm_orthogonality_U} gives
\[
    \|u_{i,\delta}\|_{L^2(U)}\longrightarrow1.
\]
Therefore
\[
    \|v_{i,\delta}\|_{L^2(U)}
    \longrightarrow1,
\]
and hence
\[
    \|u_{i,\delta}-\hat v_{i,\delta}\|_{L^2(U)}
    \longrightarrow0.
\]
It follows that, along the diagonal sequence,
\[
    u_{i,\delta_n}
    \longrightarrow u_i
    \qquad\text{in }L^2(U)
\]
for every fixed $i$.

Together with the convergence on $I\setminus U$ shown at the beginning of the proof, and the fact that both functions vanish on $\R\setminus I$, this gives
\[
    u_{i,\delta_n}
    \longrightarrow u_i
    \qquad\text{in }L^2(\R)
\]
for every $i$.

Finally, the functions $\{u_i\}_{i\geq1}$ are orthonormal, and
$u_i$ is an eigenfunction corresponding to $\lambda_i$ for every $i$. If an eigenvalue $\lambda$ has multiplicity $m$, then as the eigenvalues are listed according to multiplicity, exactly $m$ of the functions $u_i$ correspond to $\lambda$. These functions are orthonormal and therefore
form a basis of that eigenspace. Thus the family $\{u_i\}_{i\geq1}$ spans every eigenspace of the infinite-well operator. As the infinite-well operator has a complete orthonormal system of eigenfunctions in $L^2(U)$, the closed span of its eigenspaces is $L^2(U)$. It follows that $\{u_i\}_{i\geq1}$ is 
an orthonormal basis for $L^2(U)$.
\end{proof}

To complete the counterexample, let $U$ consist of three subintervals, and choose the bounded smooth potential $V$ on $U$ to be the same one as in the infinite potential well problem of Section \ref{sec_limiting_counterexample};  i.e., $V$ takes a constant value $V_i$ on $I_i$, for $i=1,2,3$, with $V_2<V_1,V_3$, so that we have a bounded triple well potential on $I$ in eigenvalue problem (\ref{evalue_problem_finite_well}). See Figure \ref{fig_finite_potential}. Recall that the first three eigenvalues are simple in this specific infinite potential well situation, as shown in Section \ref{sec_limiting_counterexample}. 
By Theorem \ref{L2_convergence_efunctions}, along a sequence
$\delta\to0^+$, an eigenfunction $u_{i,\delta}$ for eigenvalue
$\lambda_{i,\delta}$ of the finite potential well problem
(\ref{evalue_problem_finite_well}) converges in $L^2$ to an
eigenfunction $u_i$ of eigenvalue $\lambda_i$ of the infinite
potential well problem (\ref{evalue_problem_infinite_well}).

It remains to justify that this $L^2$ convergence forces two sign changes in the second eigenfunction of the finite potential well problem. After replacing $u_2$ by $-u_2$ if
necessary, choose three ordered points $x_1<x_2<x_3$ in $U$ such that
\[
    u_2(x_1)<0,\qquad u_2(x_2)>0,\qquad u_2(x_3)<0.
\]
This is possible because the infinite-well second eigenfunction constructed
in Section \ref{sec_limiting_counterexample} has two sign changes. By
continuity of $u_2$, there exist mutually disjoint compact subintervals
$J_1,J_2,J_3\subset U$, ordered from left to right, for which
\[
    \int_{J_1}u_2\,dx<0,\qquad
    \int_{J_2}u_2\,dx>0,\qquad
    \int_{J_3}u_2\,dx<0.
\]
Along the sequence $\delta\to0^+$ supplied by Theorem
\ref{L2_convergence_efunctions}, the Cauchy--Schwarz inequality gives
\[
\left|\int_{J_m}(u_{2,\delta}-u_2)\,dx\right|
\leq |J_m|^{1/2}
\|u_{2,\delta}-u_2\|_{L^2(U)}
\longrightarrow0,
\qquad m=1,2,3.
\]
Thus for all sufficiently small $\delta>0$,
\[
    \int_{J_1}u_{2,\delta}\,dx<0,\qquad
    \int_{J_2}u_{2,\delta}\,dx>0,\qquad
    \int_{J_3}u_{2,\delta}\,dx<0.
\]
This implies there are points $y_m\in J_m$, ordered $y_1<y_2<y_3$, such that
\[
    u_{2,\delta}(y_1)<0,\qquad
    u_{2,\delta}(y_2)>0,\qquad
    u_{2,\delta}(y_3)<0.
\]
As noted in the Section~\ref{ch_4_background}, the finite-well eigenfunctions are continuous, so the intermediate value theorem produces one sign change between $y_1$ and $y_2$, and another between $y_2$ and $y_3$. It follows that the finite-well second eigenfunction has at least two sign changes for a sufficiently deep finite potential well.

\section{Future directions} 
The challenge of extending the counterexample to irrational $s$ remains an open question. It would likely require a different approach to the perturbation analysis for the infinite well, since the restriction to rational $s$ arises from satisfying the analyticity condition in the Kato-Rellich perturbation theory. Whether the nodal properties (at least for the second eigenfunction) of the unperturbed nonlocal Dirichlet problem are preserved under a perturbation with a suitable convex potential is also a problem worth exploring.

\section*{Acknowledgments}
The second author (S.Chen) thanks one of their thesis examiners for suggesting the Perron-Frobenius theorem could be applied to simplify the proof of Theorem \ref{positive_evector}.

\section*{Declarations}


\begin{itemize}
\item Funding.
The results in this paper form a core part of the second author's doctoral thesis, completed at the Mathematical Sciences Institute of the Australian National University. Her work was supported by an Australian Postgraduate Award. 
\item Conflict of interest/Competing interests. The authors have no competing interests to declare that are relevant to the content of this article.
\item Ethics approval and consent to participate. Not applicable.
\item Consent for publication. Not applicable.
\item Data availability. Not applicable.
\item Materials availability. Not applicable.
\item Code availability. Not applicable.
\item Author contribution. Conceptualisation: B. Andrews; Methodology: S. Chen and B. Andrews; Investigation: S. Chen; Writing: S. Chen; Writing – Review and Editing: S. Chen and B. Andrews; Supervision: B. Andrews. All authors have read and agreed to the published version of the manuscript.
\item Declaration of generative AI and AI-assisted technologies in the manuscript preparation process. After receiving feedback from reviewers for this journal, the second author used ChatGPT5.6 Sol to carefully examine impacted arguments, and to check and edit revisions. The second author reviewed and edited the output as needed. Both authors take full responsibility for the content of the published article. 
\end{itemize}


\end{document}